\documentclass[11pt]{article}
\usepackage[affil-it]{authblk}
\usepackage[table,svgnames]{xcolor}
\usepackage{amsmath,empheq}
\usepackage{geometry}
\usepackage{graphicx}
\usepackage{comment}
\usepackage{amsfonts}
\usepackage{mathrsfs}
\usepackage{float}
\usepackage{booktabs}
\usepackage{multirow}
\usepackage{complexity}
\usepackage[ruled,vlined,noend]{algorithm2e}
\SetKwInput{KwData}{Input}
\SetKwInput{KwResult}{Output}
\usepackage{adjustbox}
\usepackage{longtable}
\usepackage{anyfontsize}
\usepackage{url}
\usepackage{hyperref}
\hypersetup{
  colorlinks   = true,
  urlcolor     = cyan,
  linkcolor    = blue,
  citecolor   = red
}
\usepackage{subcaption}
\usepackage{soul}

\newtheorem{proced}{Procedure}
\newtheorem{obs}{Observation}

\makeatletter
\def\blfootnote{\xdef\@thefnmark{}\@footnotetext}
\makeatother

\makeatletter
\def\ps@pprintTitle{%
  \let\@oddhead\@empty
  \let\@evenhead\@empty
  \let\@oddfoot\@empty
  \let\@evenfoot\@oddfoot
}
\makeatother
\title{Advanced Kernel Search approach for the MST Problem with conflicts involving affinity detection and initial solution construction\blfootnote{\textit{*Corresponding author}\\\textit{Email addresses:}\\ \texttt{fcarrabs@unisa.it} (Francesco Carrabs), \texttt{mcerulli@unisa.it} (Martina Cerulli), \texttt{dserra@unisa.it} (Domenico Serra)}}

\author[1]{Francesco Carrabs}
\author[2,*]{Martina Cerulli}
\author[1]{Domenico Serra}

\affil[1]{Department of Mathematics, University of Salerno, Fisciano, 84084, Italy}
\affil[2]{Department of Computer Science, University of Salerno, Fisciano, 84084, Italy}

\begin{document}
\maketitle
            
\begin{abstract}
The Minimum Spanning Tree Problem with Conflicts consists in finding the minimum conflict-free spanning tree of a graph, i.e., the spanning tree of minimum cost, including no pairs of edges that are in conflict. 
In this paper, we solve this problem using an enhanced Kernel Search method, which iteratively solves refined problem restrictions. Our approach addresses two central open questions in the kernel search literature: (1) how to determine the affinity between variables to ensure that the restricted problem contains variables that are as compatible as possible, meaning they are more likely to appear together in a feasible solution, and (2) how to construct an initial feasible solution quickly. To this end, we integrate the computation of independent sets from the conflict graph within the algorithm to detect affinities and effectively manage conflicts. Furthermore, we heuristically construct an initial starting point, significantly accelerating the computational process. Although our methodology is designed for MSTC, its principles could be extended to other combinatorial optimization problems with conflicts. Experimental results on benchmark instances demonstrate the efficiency and competitiveness of our approach compared to existing methods in the literature, achieving 17 new best-known values.
\end{abstract}

\section{Introduction}\label{sec:intro}
The Minimum Spanning Tree Problem (MST) consists in finding the spanning tree of an edge-weighted graph, with the minimum possible total edge weight.
In this paper, we focus on the so-called Minimum Spanning Tree Problem with Conflicts (MSTC). It is a \NP-hard variant of the classical MST, with multiple pairs of edges in conflict, where two edges are said to be in conflict if they cannot be both included in the spanning tree at the same time.
MSTC has several practical applications. 
It may be used for the design of offshore wind farm networks \cite{klein2015}, where overlapping cables (to be avoided) can be seen as conflicting edges, or in the analysis of road maps with forbidden transitions \cite{kante2013trees}, i.e., maps where some routes are not valid because, for example, at some points, it is not allowed to turn left or right. Another practical application could be the installation of an oil pipeline system that connects multiple countries, where conflicts derive from the fact that several firms and several countries are involved in the system \cite{darmann2009}. Furthermore, MSTC reveals its versatility by encompassing a well-studied problem as a special case: the Hamiltonian path problem on a directed graph \cite{gutin2006traveling}. We also present a new practical application in the context of security in distributed environments \cite{Dawood2014}. The nodes in the graph may represent the entities within the environment (users, devices, servers, or any identifiable components involved in secure interactions), and the edges may represent cryptographic connections or relationships between entities, reflecting the use of cryptographic mechanisms for secure identification, authentication, and access control. Conflicts between edges model incompatibilities or interferences in the secure communication or cryptographic processes. For instance, conflicts could arise if certain cryptographic protocols cannot be used together or if there are limitations in simultaneous access control. The MSTC in this context would be used to find the optimal set of secure communication paths or cryptographic connections that link all entities while avoiding conflicts. It ensures the most efficient and conflict-free way for entities to interact securely.

The MSTC falls within the class of optimization problems with conflict (or disjunctive) constraints such as the knapsack problem with conflicts \cite{Pferschy_2009, Coniglio2021, cerulli2020},
the bin packing problem with conflicts \cite{gendreau2004,sad2013,capua2018study,Ekici2021},
the directed rural postman problem with conflicts between nodes or arcs \cite{colombi2017directed}, the minimum matching with conflicts \cite{onc2013, Oncan2018}, the shortest path with conflicts \cite{darmann2011}, and the maximum flow problems with conflicts \cite{pfe2013, Suvak2020}.
It has been proposed in~\cite{darmann2009,darmann2011}, as ``MST with disjunctive constraints'', where such disjunctive constraints are studied in terms of the so-called \textit{conflict graph}, having nodes corresponding to the edges of the original graph, and edges representing the conflicts. The authors prove that the MSTC is strongly \NP-hard, but still polynomially solvable if every connected component of the conflict graph is a single edge (i.e., a path of length one). 
Other special cases of MSTC that can be solved in polynomial time are studied in \cite{zhang2011}, where several metaheuristic algorithms to compute approximate solutions and two exact algorithms, based on Lagrangian relaxation, are presented too. 

Other heuristics are proposed in \cite{grasp}, where a Greedy Randomized Adaptive Search Procedure (GRASP) method coupled with adaptive memory programming is presented, and in \cite{carrabs2019}, where a Multiethnic Genetic Algorithm (MEGA) is designed for the Minimum Conflict Weighted Spanning Tree variant, where the primary objective is to find a spanning tree with the minimum conflicting edge-pairs, and the second one is to minimize the weight of the tree, given it is conflict-free. Comparative studies demonstrated MEGA superiority over the tabu search algorithm proposed in \cite{zhang2011}. In \cite{cerrone2019}, another genetic algorithm is proposed, which is computationally better than both the tabu search and MEGA. Recently, in \cite{chaubey2023} an approach consisting of two phases of metaheuristic techniques is proposed for the same variant of MSTC: first, a hybrid artificial bee colony algorithm is used to detect a conflict-free spanning tree, and, second, a iterated local search that starts with this conflict-free solution is used to find the MSTC of minimum weight.

In \cite{samer2015}, a general preprocessing method and a branch and cut algorithm are proposed with the generation of cutting planes corresponding to subtour elimination constraints and odd-cycle inequalities. Another branch and cut method for the MSTC is proposed in \cite{carrabs2021bc}, where a new set of valid inequalities for the problem is introduced. To obtain upper bounds that help prune the search tree, the authors of \cite{carrabs2021bc} use the genetic algorithm proposed in \cite{carrabs2019}. 
As in \cite{zhang2011}, in \cite{gaudioso2021}, the Lagrangian relaxation of the MSTC is studied and solved with an ad hoc dual ascent procedure. In order to obtain lower bounds for MSTC, the authors of~\cite{carrabs2024} proposed a subgradient method that uses a new rule to calculate the step size based on retaining some information from the previous iteration. 

In this paper, we propose a new solution approach to the MSTC based on the so-called Kernel Search (KS) method~\cite{Angelelli2010}. It falls within the category of general-purpose methods for solving Mixed-Integer Linear Programming (MILP) problems. Initially employed in addressing challenges like the multidimensional knapsack problem \cite{Angelelli2010} and a portfolio selection problem \cite{Angelelli2012}, it has also been used to solve, among others, facility location problems with or without conflicts \cite{mansini2022,filippi2021kernel,guastaroba2012kernel}, as well as the multivehicle inventory routing problem \cite{archetti2021}. The main purpose of this algorithm is to obtain a feasible solution, hopefully of good quality, from a small set of promising variables of the MILP. The algorithm operates in two phases: in the initialization phase, the so-called kernel set is initially built using the information provided by solving the Linear Programming (LP) relaxation of the original MILP. In the subsequent improvement phase, new promising variables to add to the kernel set are iteratively identified by solving MILP subproblems of increasing size.

One of the novelties in our approach is the integration of independent sets (sets of vertices, no two of which are adjacent) from proper subgraphs of the conflict graph within the KS algorithm. This is motivated by the observation that the set of edges in any feasible solution for the MSTC corresponds to an independent set in the conflict graph. By leveraging this in both the initialization and improvement phases, we ensure that the variables included in the reduced problem exhibit high compatibility (or \textit{affinity}), thereby increasing the likelihood of forming a feasible solution. We believe this concept has the potential to be extended to other optimization problems involving conflicts where the conflict graph is utilized.
Another distinguishing aspect of our methodology is the iterative nature of the KS algorithm we devise, in which buckets are not only revisited but also merged based on independent sets identified in the conflict graph. This approach enhances the algorithm ability to manage conflicts effectively and refine the solution progressively. 
Finally, we introduce a heuristic procedure for constructing a starting point still employing independent sets of the conflict graph. This addresses one of the main challenges of the KS approach: the difficulty of finding a feasible starting solution, particularly for complex problems \cite{guastaroba2017adaptive}. These starting solutions are particularly useful for building the initial kernel set for problems like MSTC in which the optimal solution of the continuous relaxation does not provide enough information to identify promising variables.
The proposed KS algorithm is tested on benchmark instances, achieving 17 new best-known values, demonstrating its effectiveness in advancing the state of the art.

The rest of the paper is organized as follows. 
In Section~\ref{sec:MSTC}, we present the Integer Linear Programming (ILP) formulation of MSTC we consider. In Section~\ref{sec:preliminaries}, we outline the KS framework. In Section~\ref{sec:heu}, we present a heuristic finding an independent set of a given conflict graph (Section~\ref{sec:independent}) and another finding a promising set of edges in the graph (Section~\ref{sec:heuristic}). In Section~\ref{sec:kernel}, we describe the two phases of the proposed KS approach in detail. A preprocessing procedure as well as some valid inequalities for the MSTC are described, which are used to strengthen the LP relaxation of the problem. In Section~\ref{sec:tests}, we test the KS approach on the benchmark instances and compare our results with the ones obtained by other heuristics available in the literature. Finally, Section~\ref{sec:conclusion} concludes the paper.

\section{ILP Formulation of MSTC and notation}\label{sec:MSTC}
Let us consider an undirected graph $\mathcal{G}=(V,E)$, where $V=\{1,\dots,n\}$ is the set of nodes and $E=\{e_1,\dots,e_m\}$ the set of edges. The edge $e_i$ can be defined also by its two endpoints $u,v \in V$. For a given subset of nodes $S \subseteq V$, $E(S)=\{\{u,v\} \in E~|~u,v \in S\}$ is the set of edges having both the two endpoints in $S$.
Let $C$ be the set of conflicting edge pairs, which we call \textit{conflict set}: $\;C = \{(e_i, e_j)~:~e_i,e_j \in E, \text{ and } e_i \text{ is in conflict with } e_j\}.$

The general formulation for MSTC uses the following binary variables, defined for each $i$ s.t. $e_i \in E:$
\begin{equation*}\label{eq:xvar}
    		x_{i} = \begin{cases}
			1 &\text{if edge $e_i$ is included in the tree}\\
			0 &\text{otherwise.}
		\end{cases}
\end{equation*}
The ILP formulation of the MSTC, coming from the \textit{Subtour Elimination} formulation of the classical MST, reads:

\begin{subequations}\label{eq:mstc}
\begin{empheq}[left= \mbox{($\mathsf{MSTC}$)\quad } \empheqlbrace]{align}
    \min\limits_{x} & \sum_{i: e_i \in E} w_i x_i & \label{eq:mstc_obj}\\
    \text{s.t.} & \;\sum\limits_{i: e_i \in E} x_i = n - 1 & \label{eq:mstc1}\\
    & \sum\limits_{i: e_i \in E(S)} x_i \leq |S|-1 & \forall \, S \subset V, |S|\geq 3 \label{eq:mstc2}\\
    & x_i + x_j \leq 1 & \forall\, i,j:(e_i,e_j)\in C \label{eq:mstc3}\\
    & x \in \{0,1\}^{|E|}\label{eq:mstc4} &
\end{empheq}
\end{subequations}
The objective function~\eqref{eq:mstc_obj} represents the cost of the spanning tree, which we want to minimize. Constraint \eqref{eq:mstc1} imposes that $n-1$ edges must be selected since a spanning tree must include exactly $|V|-1$ edges. The exponentially many constraints~\eqref{eq:mstc2} are needed to ensure that no subtour is contained in the selected edges, i.e., no subgraph induced by $S$ contains a cycle. Finally, constraints~\eqref{eq:mstc3} impose that, for each pair of edges in conflict, at most one edge of such pair can be selected. 

In the following, we denote by LP-MSTC the LP relaxation of model~\eqref{eq:mstc}, obtained by relaxing constraints~\eqref{eq:mstc4}. In addition, we define as $\overline{\text{LP-MSTC}}$ the model LP-MSTC without the subtour elimination constraints (constraints~\eqref{eq:mstc2}).
Given a subset of edges $F \subseteq E$, we further denote by MSTC$(F)$ the ILP restriction of problem~\eqref{eq:mstc} obtained by restricting the set of variables to the edges in $F$, i.e., adding constraints 
$$x_i = 0 \quad \forall i \in E\setminus F$$ 
to formulation~\eqref{eq:mstc}. As before, LP-MSTC$(F)$ denotes the LP relaxation of MSTC$(F)$.
We further use the notation $conflicts(F)$ to represent the number of conflicts in the set of edges $F$, i.e., the cardinality of set $\{i,j\in F:(e_i,e_j)\in C \}$. 

Furthermore, to indicate the $k$-th element of a set $A$, we use the notation $A[k]$, with $A[1]$ representing the first element of $A$, and to refer to the first $k$ elements of set $A$, we use the notation $A[1:k]$. Instead, given a graph $\mathcal{H}$, a subset of its nodes $\tilde V$, and a subset of its edges $\tilde E$, we denote by $\mathcal{H}[\tilde V]$ the induced subgraph of $\mathcal{H}$ whose node set is $\tilde V$ and whose edge set consists of all of the edges of $\mathcal{H}$ that have both endpoints in $\tilde V$. Whereas, we denote by $\mathcal{H}[\tilde E]$ the induced subgraph of $\mathcal{H}$ whose edge set is $\tilde E$ and whose node set consists of all of the endpoints of $\tilde E$ in $\mathcal{H}$. 

To conclude this section, we recall some important definitions. Let $N(u)=\{v\in V : \{u,v\}\in E\}$ be the set of \textit{neighbors} of node $u$ that is the set of nodes adjacent to $u$ in $\mathcal{G}$. An \textit{independent set} is a set of nodes in a graph, no two of which are adjacent. 
Given the graph $\mathcal{G}=(V, E)$, its \textit{conflict graph} $\hat{\mathcal{G}}= (\hat{V}, \hat{E})$ is a graph which has a node for each edge in $\mathcal{G}$, and where we represent each conﬂict constraint by an edge connecting the corresponding vertices in $\hat{\mathcal{G}}$. Note that each conﬂict-free spanning tree in $\mathcal{G}$ is a subset of $\hat{V}$ which corresponds to an independent set of $\hat{\mathcal{G}}$.

\section{Preliminaries}\label{sec:preliminaries}
In this section, we briefly outline the KS framework used to solve the MSTC problem and introduce two algorithms employed within it. 
Given the mathematical formulation~\eqref{eq:mstc} of MSTC, let $\Omega$ be the set of \textit{ideal} variables, i.e., the variables taking a value equal to one in the optimal solution of MSTC. It is easy to see that the resolution of the subproblem MSTC$(\Omega)$ provides the optimal solution to the MSTC problem: all the $|\Omega| = n - 1$ variables must be set to 1. 
As $\Omega$ is not known, the KS algorithm consists in identifying a set {$\Lambda$} of \textit{promising} variables of this model, i.e., variables that have a high probability of being non-zero in the optimal solution of the problem. This is because, if $\Omega \subseteq \Lambda$, the resolution of MSTC$(\Lambda)$ provides the optimal solution to the MSTC problem. The starting set of promising variables is used to form the \textit{kernel set}, while the remaining ones are partitioned into several subsets called \textit{buckets}. We define $\Lambda^0$ as the set of variables indices in the initial kernel set, and $B_k$ as the set of variables indices in the $k$-th bucket, for $k=1,\dots,b$. 

The partitioning of the indices of the variables into the kernel set and the buckets is carried out by KS according to the optimal solution of LP-MSTC, (eventually) a feasible MSTC solution, as well as an independent set of a proper subgraph of the conflict graph $\hat{\mathcal{G}}$. After the construction of the kernel set and the buckets, KS seeks an initial feasible solution by solving the MSTC restricted to the variables having indices inside the kernel set only, i.e., solving MSTC($\Lambda^0$). For this reason, the initial kernel size $|\Lambda^0|$ should be, at the same time, large enough to contain as many variables as possible among the ones that are likely to be part of an optimal solution, and not too large to ensure that MSTC($\Lambda^0$) can be efficiently solved to optimality.

The kernel set is then updated throughout the KS algorithm ($\Lambda^k$ is the kernel set at the end of iteration $k$), by taking into account the buckets. At each iteration, indeed, the MILP formulation, restricted to the current kernel set and one or multiple buckets from the sequence, is solved, trying to obtain solutions of increasing quality and thus new variables to insert into the kernel set.

To sum up, the KS algorithm is divided into two phases:
\begin{itemize}
     \item \textbf {Initialization phase}: The initial set $\Lambda^0$ and the bucket sequence $\mathcal{B} = \{B_1, \dots B_{b}\}$ are constructed. 
     Then, the MSTC($\Lambda^0$) is solved and, if a feasible solution is found, its objective function value is set as the current upper bound $w^*$ of the optimal solution we are looking for. 
     \item \textbf{Improvement phase}: At each iteration $k$, the MSTC($\Lambda^{k-1}\cup B_k$) is solved. If this subproblem is feasible, its optimal solution $x^k$ and relative value $w^k$ are used to, eventually, update the incumbent solution and its value, respectively. Moreover, all the variables of $B_k$, greater than zero in $x^k$, are moved into the kernel set, generating the new kernel set $\Lambda^k$. The iteration over the buckets may be performed more than once, considering new buckets obtained from the union of the original buckets.
     \end{itemize} 

All the parameters that will be involved in the KS algorithm are listed in Table~\ref{tab:parameters}.
\begin{table}[ht]
\centering  
\resizebox*{\textwidth}{!}{
\begin{tabular}{|rl|}
\hline
\rowcolor[HTML]{C0C0C0} 
\textbf{Parameter} & \textbf{Description} \\
$\bar h$ & maximum number of for-loop iterations of the \texttt{StartingSolution} algorithm\\
\rowcolor{gray!10} 
$\bar t$ & maximum number of while-loop iterations of the \texttt{StartingSolution} algorithm \\
$b$ & number of buckets \\
\rowcolor{gray!10} 
$d$ & cardinality of each bucket \\
$K$ & cardinality of the set which is used to create the first kernel set\\
\rowcolor{gray!10} 
$P$ & number of iterations performed over the buckets in the improvement phase \\
\multirow{2}{*}{$p \in \{1, P\}$} & index of the iteration performed over the buckets in the improvement phase representing the type of bucket \\
& added to the kernel set (outer iterations) \\
\rowcolor{gray!10} 
$\Delta$ & maximum number of buckets added to the kernel at each iteration, without finding a feasible MSTC solution\\ 
$\bar k \in \{1,\Delta\}$ & index representing the number of buckets already added to the kernel set without finding a feasible MSTC solution\\
\rowcolor{gray!10} 
$k \in \{1,b\}$ & index representing the bucket added to the kernel set (inner iterations) \\
$\Lambda^k$ & kernel set at iteration $k$ \\
\rowcolor{gray!10} 
$B_k$ & $k$-th bucket with $k \in \{1,b\}$\\
$\mathcal{B}$ & set of buckets \\
\hline
    \end{tabular}}
\caption{List of the parameters involved in the KS algorithm.}
\label{tab:parameters}
\end{table}

\section{Embedded heuristics}\label{sec:heu}
The effectiveness of the KS approach depends on the quality of the initial kernel set. This set is usually generated using the information retrieved from the LP relaxation of the model. However, for some problems, this information alone is not sufficient to ensure a high-quality initial kernel set. For this reason, an alternative approach, recently investigated in \cite{carrabs2025}, consists of applying fast heuristic methods to generate feasible solutions to the problem. These solutions provide additional information, which, together with the LP relaxation solutions, can be used to construct the initial kernel set. To this end, in this section, we introduce two heuristics.

The first heuristic identifies independent sets in the conflict graph, ensuring that the selected edges do not induce cycles in the original graph. While this guarantees a conflict-free selection, the resulting subgraph may not be connected.

To improve connectivity, the second heuristic uses the first one but iteratively extends its outputs by incorporating additional edges. It can yield different outcomes: a conflict-free spanning forest, a spanning tree that includes some conflicting edges, or an entirely conflict-free spanning tree. It leverages the first heuristic within its iterations to maintain as many conflict-free edges as possible. Through this approach, the second heuristic effectively attempts to produce a feasible MSTC solution.

\subsection{Computing an independent set of the conflict graph}\label{sec:independent}
In this section, we describe a greedy algorithm, \texttt{IndependentSet}, whose aim is to compute an independent set $I$ of the conflict graph $\mathcal{\hat{G}}$ that does not create cycles inside the original graph $\mathcal{G}$. The usefulness of this algorithm derives from the following observation.

\begin{obs}
    Any feasible MSTC solution is composed of a set of $n-1$ edges that: (i) does not create cycles in $\mathcal{G}$; (ii) corresponds to an independent set of $\mathcal{\hat{G}}$.
\end{obs} 

Our idea is to use the \texttt{IndependentSet} algorithm during both phases of KS to find promising variables. It is worth noting that the idea of using the maximum independent set can be extended to any optimization problem with conflicts because its feasible solutions must be independent sets of the conflict graph.

The pseudocode of the \texttt{IndependentSet} algorithm is shown in Algorithm~\ref{algo:indset}. The algorithm takes in input the graphs $\mathcal{G}$ of MSTC and a subgraph $\mathcal{H}$ of the conflict graph $\mathcal{\hat{G}}$. The independent set $I$ is initialized to the empty set on line~\ref{initialization-indset} together with sets $\bar{V}$ and $\bar{E}$. The set of vertices $W$ is initialized with all the vertices of $\mathcal{H}$ (line~\ref{initialization-W}), and the while loop (lines \ref{while-is-start}-\ref{while-is-end}) is repeated until $W$ is not empty.

\noindent\begin{minipage}{0.95\textwidth}%
{\small\begin{algorithm}[H]\caption{\texttt{IndependentSet}}
    \LinesNumbered
    \SetAlgoLined 
    \label{algo:indset}
    \KwData{Graph $\mathcal{G} = (V,E)$, graph $\mathcal{H}=(V_H, E_H)$.} 
    \KwResult{Independent set $I$ of $\mathcal{H}$.} 
    $I \gets \emptyset, \bar{V} \gets \emptyset, \bar{E} \gets \emptyset$.\label{initialization-indset}\\
    $W \gets V_H.$\label{initialization-W}\\
    \While{$W \neq \emptyset$\label{while-is-start}}{
    Select a node $i \in W$ with minimum degree in $\mathcal{H}[W]$, and let $e_i = \{u, v\} \in E$ be the edge of $\mathcal{G}$ associated to $i$. \label{selectnode}\\
    \eIf{ $\mathcal{G}[\bar{E}\cup\{e_i\}]$ does not contain cycles\label{check-cycles-start}}
    {$I \gets I \cup \{i\}.$\label{update1} \\
    $\bar{V} \gets \bar{V} \cup \{u,v\}$, $\bar{E} \gets  \bar{E} \cup \{e_i\}.$\label{update2}\\
    $W_i \gets $ set of neighbors of $i$ in $\mathcal{H}$, including $i$. 
    \\$W \gets W \setminus W_i.$\label{update3}
    }
    {
    $W \gets W \setminus \{i\}.$
    }\label{check-cycles-end}
    }\label{while-is-end}
    \Return{$I$}
\end{algorithm}}
\end{minipage}

The algorithm selects a node $e_i$ of $W$ having the minimum degree in the induced subgraph $\mathcal{H}[W]$ (line~\ref{selectnode}). Let us suppose that $e_i$ corresponds to the edge $i$ in $\mathcal{G}$ whose extremes are $u$ and $v$. The \texttt{IndependentSet} algorithm checks if the subgraph of $\mathcal{G}$ induced by $\bar{E}\cup\{e_i\}$ contains cycles (line~\ref{check-cycles-start}). If this is the case, the node $i$ is removed from W (line~\ref{check-cycles-end}). Otherwise, $i$ is added to the independent set $I$ (line~\ref{update1}), and $\bar{E}$ and $\bar{V}$ are updated by adding the edge $e_i$ and its extremes $u$ and $v$, respectively (line~\ref{update2}). Finally, $i$ and its neighbors are removed from $W$ (line~\ref{update3}).
Once $W$ becomes empty, the algorithm returns the built independent set $I$.

Since the \texttt{IndependentSet} algorithm selects edges that are not in conflict with each other and that do not generate cycles in the starting graph $\mathcal{G}$, then these edges induce either a feasible solution of the MSTC or a set of conflict-free \textit{subtrees} of $\mathcal{G}$.
This means that, by selecting the variables according to the solution of the \texttt{IndependentSet} algorithm, we increase the probability of including edges in the kernel set that can generate feasible MSTC solutions. We further use the \texttt{IndependentSet} algorithm within the improvement phase to decide how to merge the buckets if multiple of them are supposed to be considered at a time.

\subsection{Finding a feasible starting point}\label{sec:heuristic}
In this section, we present the \texttt{StartingSolution} heuristic procedure that tries to generate a feasible MSTC solution, leveraging the set of edges returned by the \texttt{IndependentSet} algorithm, which is conflict-free, but not necessarily connected. The \texttt{StartingSolution} pseudo-code is shown in Algorithm~\ref{algo:heuristic}.

\noindent\begin{minipage}{0.95\textwidth}%
{\small\begin{algorithm}[H]\caption{\texttt{StartingSolution}}
    \LinesNumbered
    \SetAlgoLined 
    \label{algo:heuristic}
    \KwData{Graph $\mathcal{G} = (V,E)$, conflict graph $\hat{\mathcal{G}}$, initial edge set $S^0$, two integer parameters $\bar{h},\bar{t}>1$.} 
    \KwResult{Either a conflict-free spanning tree or an empty set.} 
    $S \gets S^0$. \\
    $\bar E \gets \emptyset$.\\
    $\bar w\gets +\infty$.\\
    \For{$h = 1,\dots \bar h$ \label{for-start}}{
    Launch the Kruskal algorithm on $\mathcal{G}[S]$ obtaining a set of edges $E'$.\label{kruskal}\\
    $\hat E \gets\{e_i \in E': \exists\; e_j \in E' \; s.t.\; (e_i,e_j) \in C\}$.\label{hatE}\\
    $t \gets 0$ \\
    \While{$|\hat E| \neq 0$ and $t\leq \bar t$ \label{while-start}}{
    Compute the independent set $I'$ over graph $\hat{\mathcal{G}}[\hat E]$ through the \texttt{IndependentSet} algorithm. \label{algo1-1}\\
    $E' \gets E' \setminus (\hat E \setminus I')$. \label{updateE}\\ 
    Compute the independent set $I$ over graph $\hat{\mathcal{G}}$ through the \texttt{IndependentSet} algorithm, imposing that $E' \subseteq I$. \label{algo1-2}\\
    \eIf{$I$ is a spanning tree}{$\hat{E}\gets\emptyset$ and $E' \gets I$\\
    \textbf{break} \label{break}}
    {Assign a random weight ${w'}_e$ to each edge $E$. \label{graphaugmentation-start} \\
    Connect $I$ with the edges in set $\bar I \subseteq E\setminus I$ such that $\sum_{e \in \bar I} {w'}_e$ is minimized, obtaining the tree $T$. \label{graphaugmentation-end}\\
    $E' \gets T.$\\
    $\hat E =\{e_i \in T: \exists\; e_j \in T \; s.t.\; (e_i,e_j) \in C\}$.\label{stepconf} \\}
    $t \gets t+1$ \\
    }
    $S \gets S \cup E'.$ \label{updateS}\\
    \If{$\hat{E} = \emptyset$ and $\sum\limits_{e \in E'} w_e < \bar w$\label{if-start}}{$\bar E \gets E'$\\ $\bar w \gets \sum\limits_{e \in E'} w_e$ \\ \lIf{$t = 0$}{\textbf{break}}\label{if-end}}
    }
    \Return{$S,\bar E$}
\end{algorithm}}
\end{minipage}

The procedure makes use of Kruskal's algorithm \cite{kruskal1956shortest}. The process starts with a set of edges $S=S^0$. A for-loop of maximum $\bar h$ iterations begins, and Kruskal's algorithm is applied to the induced subgraph $\mathcal{G}[S]$ (line~\ref{for-start}-\ref{kruskal}). If the resulting forest $E'$ (or tree if $\mathcal{G}[S]$ is connected) contains conflicting edges, a while-loop of at most $\bar t$ iterations is initiated (line~\ref{while-start}). Let $\hat{E}$ be the edges of $E'$ that are in conflict with each other. Within each iteration of this while-loop, an independent set $I'$ is computed over the conflict graph obtained from the conflicting edges in $E'$, using the \texttt{IndependentSet} algorithm (line~\ref{algo1-1}). The edges in $\hat{E}$ that are not part of $I'$ are removed from the forest $E'$ (line~\ref{updateE}). 
A new independent set $I$, containing $E'$, is computed by the \texttt{IndependentSet} algorithm on the whole conflict graph $\hat{G}$ (line~\ref{algo1-2}).
If $I$ forms a spanning tree, a feasible solution to the MSTC problem is obtained, and the while-loop terminates (line~\ref{break}). Otherwise, additional edges are added to complete the tree through a graph augmentation process (lines~\ref{graphaugmentation-start}--\ref{graphaugmentation-end}). This graph augmentation procedure finds a minimum spanning tree $T$, containing $I$, according to randomly assigned weights. 
If any conflicting edge pair is detected in $T$ (line~\ref{stepconf}), and the maximum number of iterations for the while-loop has not been reached, a new iteration begins. Otherwise, the conflict-free tree is added to the set $S$ (line~\ref{updateS}) and eventually used to update the incumbent solution (lines~\ref{if-start}--\ref{if-end}). With this enlarged set $S$, the for-loop then proceeds to the next iteration, unless either the maximum number of iterations has been reached or $t=0$ (line~\ref{if-end}). This last condition holds when $\hat{E}=\emptyset$ on line~\ref{hatE} and, thus, the while-loop of line~\ref{while-start} is not performed. Consequently, $E'\subseteq S$ and therefore $S$ (line~\ref{updateS}) could not change for the remaining iterations of the for-loop. 
In the end, the heuristic will return as $\bar E$ one of the following outcomes: either an empty set or a conflict-free spanning tree. 

As the heuristic will be used in the initialization phase of the KS algorithm to create the initial kernel set $\Lambda^0,$ and the returned set of edges may be also empty, the heuristic further returns the set $S$ of analyzed edges.

\section{Kernel Search Approach}\label{sec:kernel}
We are now ready to describe in detail our KS approach. In Section~\ref{sec:initialization}, we focus on the initialization phase, which aims at defining the starting kernel set and the starting set of buckets. In Section~\ref{sec:extension}, we describe the improvement phase of the algorithm, in which better and better restrictions of the MSTC problem are solved. In the rest of the paper, we will refer to the entire KS algorithm, including both the initialization and extension phases, as \texttt{KS}.

\subsection{Initialization phase}\label{sec:initialization}
The construction of the starting kernel set $\Lambda^0$ and of the set of buckets $\mathcal{B}$ is carried out by Algorithm~\ref{algo:kernel_search-if}. 
The algorithm takes as input two positive parameters $K$ and $d$ (see Table~\ref{tab:parameters}) and returns both the sets $\Lambda^0$, and $\mathcal{B}$, as well as a starting solution $x^0$ together with its value $w^0.$ 
{\small\begin{algorithm}[!ht]\caption{\texttt{KS: Initialization Phase}}
    \LinesNumbered
    \SetAlgoLined 
    \label{algo:kernel_search-if}
    \KwData{Graph $\mathcal{G}=(V,E),$ conflict graph $\hat{\mathcal{G}}$, two scalars $K>0, d>0$.}
    \KwResult{Initial kernel set $\Lambda^0$, initial bucket set $\mathcal{B}$, the first MSTC solution found $x^0$, and its value $w^0.$} 
    
    Compute the optimal solution $\bar{x}'$ of the LP-MSTC relaxation. \label{compute_root_LP}\\
    Compute the optimal solution $\bar{x}''$ of the $\overline{\text{LP-MSTC}}$ relaxation. \label{compute_LP}\\
    \If{conflicts$(\{i:\bar{x}_{i}' > 0\})$ $<$ conflicts$(\{i:\bar{x}_{i}'' > 0\})$}{\label{lp_if_start}
    $\bar{x} \gets \bar{x}'$
    }
    \Else{
        $\bar{x} \gets \bar{x}''$
    }\label{lp_if_end}
    
    $N \gets\{i \in E : \bar x_i>0\}$. \label{lpmstc} \\                            
    Sort indices $i \in N$ in non-increasing order of $\bar x_i$ values. \label{sortlp}\\
    $K \gets \min \{K, |N|\}$. \label{defineK}\\ 
    $N_K \gets N[1:K]$. \label{defineconfl} \\ 
    Compute the sets $\bar{E}$ and $S$ through the \texttt{StartingSolution} algorithm, with starting edge set $S^0 := N_K$.\label{compE}\\
    Compute the independent set $I$ over the graph $\mathcal{H} :=\hat{\mathcal{G}}[S]$, through the \texttt{IndependentSet} algorithm.\label{compI}\\ 
    $\Lambda^0 \gets I \cup \bar{E}$, and $L \gets E \setminus \Lambda^0$. \label{lambda0}\\
    Sort indices $i \in L$ according to Procedure~\ref{ordering-proced}. \label{sortL}\\
    \If{$|\Lambda^0| < K$\label{finishkernel-start}}{$\gamma = K-|\Lambda^0|$\\
    $\bar L = L[1:\gamma]$\\
    $\Lambda^0 \gets\Lambda^0 \cup \bar L$\\
    $L \gets L \setminus \bar L$\\}\label{finishkernel-end}
    $k \gets 1$\\
    \While{$L\neq \emptyset$ \label{populateB-start}}{
    $d \gets \min\{d,|L|\}$\\
    $B_k \gets L[1:d]$ \\
    $\mathcal{B} \gets \mathcal{B} \cup B_k$\\
    $L \gets L \setminus B_k$\\
    $k \gets k + 1$}\label{populateB-end}
    Solve MSTC$(\Lambda^0)$.\label{x0-start}\\
    \eIf{MSTC$(\Lambda^0)$ is feasible}{Let $x^0$ be the optimal solution of MSTC$(\Lambda^0)$, and $w^0$ its optimal value.}
    {$x^0 \gets 0$, and 
    $w^0 \gets +\infty$. \label{mstc-infeas}}\label{x0-end}
    \Return{$\Lambda^0, \mathcal{B}, x^0, w^0.$}
\end{algorithm}}

First, the optimal solutions of both the LP-MSTC relaxation and the $\overline{\text{LP-MSTC}}$ relaxation are computed (lines~\ref{compute_root_LP}--\ref{compute_LP}). Then, the relaxed solution that results in the greatest number of compatible variables (i.e., those that are more likely to form an independent set in the conflict graph) is selected using the condition in lines~\ref{lp_if_start}--\ref{lp_if_end}. Specifically, the solution with fewer conflicts among variables with positive values is chosen and denoted as $\bar x$. Next, the set $N$ is constructed by including the indices of the variables with positive values in $\bar x$ (line~\ref{lpmstc}), and  $K$ is set in such a way that it is not greater than the cardinality of $N$ (line~\ref{defineK})

The subset of the first $K$ elements of $N$ is defined as $N_K$ (line~\ref{defineconfl}). 
In line~\ref{compE}, the \texttt{StartingSolution} algorithm proposed in Algorithm~\ref{algo:heuristic} is executed with $N_K$ as starting edge set $S^0.$
The heuristic returns two edge sets: $\bar E$ and $S$. The set $\bar{E}$ is either a conflict-free spanning tree, or an empty set. The set $S$ is a superset of the set $N_K$, which is used to define the graph $\mathcal{H}$ needed as input to the \texttt{IndependentSet} algorithm. Indeed, in line~\ref{compI}, a feasible solution $I$ of the maximum independent set problem on the conflict subgraph $\hat{\mathcal{G}}[S]$ is computed through the greedy procedure presented in Algorithm~\ref{algo:indset}, by setting $\mathcal{H}=\hat{\mathcal{G}}[S]$. This procedure guarantees that the returned set $I$ does not generate cycles in the starting graph $\mathcal{G}$.

Given the two edge sets $I$ and $\bar{E}$, coming from the \texttt{IndependentSet} algorithm and the \texttt{StartingSolution} algorithm, respectively, in line~\ref{lambda0} the initial kernel set $\Lambda^0$ is populated with the indices of the variables associated with the edges in $I \cup \bar{E}$, and the remaining indices are used to define set $L$, which is used first to complete the set $\Lambda^0$ if its cardinality is smaller than $K$, and then to populate the buckets.
The elements of set $L$ are sorted in line~\ref{sortL} according to the following procedure:
\begin{proced}\label{ordering-proced}
The elements of $L$ are sorted sequentially according to the following criteria: $(i)$ non-decreasing order of the number of conflicts with the edges in $\Lambda^0$; $(ii)$ non-increasing order of value in the solution of LP-MSTC; $(iii)$ non-increasing order of the reduced cost coefficients. If two variables are tied according to a given criterion, the next criterion is used to determine priority. If there is still a tie after all three criteria are applied, the variable with the larger index $i$ is given priority.
\end{proced}
In lines~\ref{finishkernel-start}--\ref{finishkernel-end}, the initial kernel set is eventually completed with $K-|\Lambda^0|$ variables in $L$. Then, in lines~\ref{populateB-start}--\ref{populateB-end}, the remaining variables in $L$ are equally distributed to populate the buckets $B_1, \dots, B_b,$ each of them of cardinality $d$ (except for the last bucket, the size of which could be smaller). The set $\mathcal{B}$ is defined as the set of buckets, and $b = |\mathcal{B}|$. In lines~\ref{x0-start}--\ref{x0-end}, the initialization phase concludes by solving the MSTC($\Lambda^0$) formulation. If it is feasible, the algorithm sets $x^0$ and $w^0$ to its optimal solution and value, respectively. Otherwise, it sets $x^0 = 0$ and $w^0 = +\infty$. Finally, the algorithm returns $\Lambda^0$, $\mathcal{B}$, $x^0$, and $w^0$.

\subsection{Improvement phase}\label{sec:extension}
With the initial kernel set $\Lambda^0$, the initial set of buckets $\mathcal{B}$, as well as the starting solution $x^0$ (which is the initial best-known solution $x^*$), and its value $w^0$ (which is the best known upper bound $w^*$), returned by the initialization phase (Algorithm~\ref{algo:kernel_search-if}), the improvement phase, described in Algorithm~\ref{algo:kernel_search-ef}, starts. We recall that, if the initialization phase does not provide any feasible solution, then $x^0$ is set to $0$ and $w^0$ is set to $+\infty$ (see line~\ref{mstc-infeas} of Algorithm~\ref{algo:kernel_search-if}).
{\small\begin{algorithm}[!ht]
\caption{\texttt{KS: Improvement Phase}}
    \LinesNumbered
    \SetAlgoLined 
    \label{algo:kernel_search-ef}
    \KwData{The solution $x^0$, its value $w^0,$ sets $\Lambda^0, \mathcal{B} =\{B_1,\dots,B_b\},$ and scalars $\Delta >0,P>0.$ \\
    \KwResult{Best solution found $x^*$ and its value $w^*$.}
    $\tilde{\mathcal{B}} \gets \mathcal{B}, x^* \gets x^0, w^* \gets w^0$, $p \gets 1$, $\bar k \gets 0$.\\
    \While{$p \leq P$\label{whilep-start}}{
        \If{$p > 1$\label{updatebuckets-start}}{$\tilde{\mathcal{B}} \gets $ \texttt{EnlargedBuckets}$(\tilde{B},\Lambda^0).\quad$ // Update the set of enlarged buckets \label{enlarging}}\label{updatebuckets-end} 
        $k\gets1.$\\
        \While{$k \leq b$\label{whilek-start}}{
            \eIf{MSTC$(\Lambda^{k-1} \cup \tilde B_k)$ is feasible}
                {Compute the optimal solution $x^k$ of MSTC$(\Lambda^{k-1} \cup \tilde B_k)$, and its optimal value $w^k$.\label{solverestr}\\
                $x^* \gets x^k$, $w^* \gets w^k$.\label{updateincumbent}\\
                $\tilde{B}_k^+\gets \{i \in \tilde B_k : x^k_i >0 \}$.\label{positivevar}\\
                $\Lambda^{k} \gets \Lambda^{k-1} \cup \tilde{B}_k^+$.\label{updatekernel}\\ 
                $\tilde{B}_k \gets \tilde{B}_k \setminus \tilde{B}_k^+$.\label{updatebucket}\\
                $\bar k  \gets 0$.\label{newincumbent}}
                {$\bar{k}  \gets \bar{k} + 1$.\label{noincumbent}}
            $k \gets k + 1$.\\
            \If{$\bar{k}\geq \Delta$ and $w^* < +\infty$\label{stopping-start}}{ $k\gets b +1$, $p\gets P+1$.}\label{stopping-end}
        }\label{whilek-end}
        $\Lambda^0 \gets \Lambda^k$.\label{updatekernel0}\\
        $p \gets p +1$.}\label{whilep-end}
    \Return $w^*,x^*$\label{ext-return}
    }
\end{algorithm}}

During the improvement phase, both the kernel set and the set of buckets are updated. The \textit{outer while loop} of lines~\ref{whilep-start}--\ref{whilep-end} manages the number of runs that the algorithm carries out on the buckets and, according to the value of the index $p$, how the buckets are combined to create the restricted problems (lines~\ref{updatebuckets-start}--\ref{updatebuckets-end}). In particular, when $p=1$, the restricted problems are generated by considering one bucket from the original set of buckets $\mathcal{B}$ at a time; if $p > 1$, multiple buckets, properly selected, are considered at a time. 
The list of \textit{enlarged} buckets $\tilde{\mathcal{B}}$ is updated, according to the $p$ value, by the \texttt{EnlargedBuckets} algorithm described in Section~\ref{sec:enlarging-buckets}. 

At each iteration $p$ of the outer while loop, starting from an initial kernel set $\Lambda^0$, $b = |\tilde{\mathcal{B}}|$ restrictions of the MSCT problem are solved, according to the \textit{inner while loop} of lines~\ref{whilek-start}--\ref{whilek-end}. At each iteration $k$ of this while loop, following the KS framework, the kernel set is, eventually, updated by embedding in it variables coming from the $k$-th bucket. In particular, at iteration $k$, the algorithm solves the restricted problem $MSTC(\Lambda^{k-1}\cup \tilde{B}_k)$ (line~\ref{solverestr}) whose variables are the ones inside the current kernel set $\Lambda^{k-1}$ and the current bucket $\tilde{B}_k$ (coming from the \texttt{EnlargedBuckets} algorithm). This restricted problem is solved with the additional constraints~$\sum_{i \in \tilde B_k}{x_{i}} \geq 1$ and $\sum_{i: e_i \in E}w_i x_i \leq w^*$ to ensure that at least one new variable from the current bucket is used and that the value of the new solution $x^k$ is at least as good as the one of the current incumbent $x^*$, respectively.

If $MSTC(\Lambda^{k-1}\cup \tilde{B}_k)$ is feasible (i.e., its optimal value $w^k$ is better than or equal to the current upper bound $w^*$), then the incumbent solution $x^*$ and its value $w^*$ are updated accordingly (line~\ref{updateincumbent}).
The subset $\tilde{B}_k^+$ of variables of the current bucket $\tilde B_k$ which are nonzero in the solution $x^k$, defined in line~\ref{positivevar}, is added to the kernel set (line~\ref{updatekernel}) and removed from the corresponding bucket (line~\ref{updatebucket}), i.e., the indices belonging to subset $\tilde{B}_k^+$ are added to $\Lambda^{k-1}$, and removed from $\tilde{B}_k$. Finally, the parameter $\bar{k}$, which is needed for the stopping criterion of the algorithm, is set to zero, because a new incumbent solution was found (line~\ref{newincumbent}).
If, instead, the restriction is infeasible, i.e., $w^k \geq w^*$, we increase parameter $\bar k$ (line~\ref{noincumbent}). Indeed, if $\bar k \geq \Delta,$ i.e., for at least $\Delta$ steps, the solution has not been improved, and at least a feasible solution has been found so far, the algorithm updates the parameters $k$ and $p$ (lines~\ref{stopping-start}--\ref{stopping-end}) to terminate the two while loops. Otherwise, the kernel set $\Lambda^0$ is initialized to the $k$-th kernel set (line~\ref{updatekernel0}), and a new iteration of the outer while loop starts.

In the last line, the algorithm returns the best solution $x^*$ found and its objective function value $w^*$ (line~\ref{ext-return}).

\subsubsection{Enlarging the buckets}\label{sec:enlarging-buckets}
Let us describe the \texttt{EnlargedBuckets} algorithm, whose pseudo-code is reported in Algorithm~\ref{algo:enlarged-bucket}, which is the method called at the beginning of each iteration of the outer while loop, at line~\ref{enlarging}, of Algorithm~\ref{algo:kernel_search-ef}. 
It starts by initializing the set $\hat{\mathcal{B}}$, used to define the new bucket set, to the empty set (line~\ref{init-B}). Two nested for loops (lines~\ref{for1-start}--\ref{for1-end}) are used to select all the possible couples $A_{lt} = \tilde B_l \cup \tilde B_t$ of buckets from $\tilde{\mathcal{B}}$. For each couple $A_{lt}$, the algorithm builds the conflict subgraph $\mathcal{H}$ whose nodes $V_H$ are all the edges in $\Lambda^0\cup A_{lt}$ (line~\ref{defineconfl-lt}). At line~\ref{computei}, an independent set of $\mathcal{H}$ is computed by using the \texttt{IndependentSet} algorithm. The cardinality of this independent set gives us an estimation of the compatibility between the edges inside the set $\Lambda^0 \cup A_{lt}$. Indeed, the greater this cardinality is, the higher it will be the probability of producing feasible solutions of MSTC by solving the restricted problem $MSTC(\Lambda^0\cup A_{lt})$. For this reason, given the set $\hat{A}$ of all the bucket pairs $A_{lt}$, these pairs are sorted, in non-increasing order of the cardinality of the associated sets $I_{lt}$ (lines~\ref{defineA}-\ref{orderA}). At each iteration $k$ of the while loop at lines~\ref{whileA-start}--\ref{whileA-end} (with $k$ initialized to $1$ in line~\ref{init-k}), the $k$-th couple of $\hat{A}$, say $A_{lt}=\tilde{B_l}\cup \tilde{B_t}$, is considered. The buckets $\tilde{B_l}$ and $\tilde{B_t}$ are removed from the set $\tilde{\mathcal{B}}$ and all the couples containing either $\tilde{B_l}$ or $\tilde{B_t}$ are removed from $\hat{A}$ (line~\ref{operations-A}). Furthermore, the new enlarged bucket $\tilde{B_l}\cup \tilde{B_t}$ is added to $\hat{\mathcal{B}}$ (line~\ref{operations-B}). If the number of current buckets is odd, the last bucket is added to $\hat{\mathcal{B}}$ (lines~\ref{Bodd-start}--\ref{Bodd-end}), that is, it will be again considered alone in Algorithm~\ref{algo:kernel_search-ef}. In the last two lines, the algorithm assigns $\hat{\mathcal{B}}$ to $\tilde{\mathcal{B}}$ and returns this last set of enlarged buckets. 

{\small\begin{algorithm}[!ht]\caption{\texttt{EnlargedBuckets}}
    \LinesNumbered
    \SetAlgoLined 
    \label{algo:enlarged-bucket}
    \KwData{Sets $\tilde{\mathcal{B}}=\{\tilde{B}_1,\dots,\tilde{B}_b\}$, with $b = |\tilde{\mathcal{B}}|$, $\Lambda^0$.}
    \KwResult{Updated set of enlarged buckets $\tilde{\mathcal{B}} = \{\tilde{B}_1,\dots,\tilde{B}_b\}$.} 
    $\hat{\mathcal{B}} \gets \emptyset$. \label{init-B}\\
    \For{$l = 1 \dots b-1$\label{for1-start}}{
            \For{$t = l+1\dots b$}{
            $A_{lt} \gets \tilde{B}_l \cup \tilde{B}_t$.\\
            $V_H\gets\Lambda^0 \cup A_{lt}$, and $E_H \gets \{ (i,j)~:~i,j \in V_H , \text{ and } (e_i , e_j ) \in C \}$.\label{defineconfl-lt}\\
            Compute the independent set $I_{lt}$ over graph $\mathcal{H}=(V_H,E_H)$ through the \texttt{IndependentSet} algorithm.\label{computei}\\
            }
            }\label{for1-end}
        Define the set $\hat A$ as the set of all the obtained bucket pairs $A_{lt}$. \label{defineA} \\
        Order set $\hat A$ elements in non-increasing order of cardinality of the associated sets~$I_{lt}$. \label{orderA}\\
        $k\gets1$.\label{init-k}\\
        \While{$\hat A \neq \emptyset$\label{whileA-start}}{
            Being $\hat A[k] = A_{lt} $, $\hat{A} \gets\hat{A} \setminus \{A_{ij}\} \, \forall\, i,j ~:~ i=l \vee j=t$, and $\tilde{\mathcal{B}} \gets \tilde{\mathcal{B}} \setminus \{\tilde{B}_l \cup \tilde{B}_t\}$. \label{operations-A}\\
            $\hat{\mathcal{B}} \gets \hat{\mathcal{B}} \cup \hat A[k]$.\label{operations-B}\\
            $k \gets k + 1$. \\
        }\label{whileA-end}
        \If{$\tilde{\mathcal{B}} \neq \emptyset$\label{Bodd-start}}
        {$\hat{\mathcal{B}} \gets \hat{\mathcal{B}} \cup \tilde{\mathcal{B}}$.\\ 
        }  \label{Bodd-end}
  $\tilde{\mathcal{B}}\gets\hat{\mathcal{B}}.$\\
    \Return{$\tilde{\mathcal{B}}$}
\end{algorithm}}

\subsection{Valid inequalities and preprocessing procedure}
In this section, first, we describe some valid inequalities for the MSTC, taken from \cite{carrabs2021bc}, that we add to our model in order to speed up the resolution of the restricted problems solved throughout the KS algorithm. Second, we present a preprocessing procedure proposed in \cite{samer2015} which sets some variables $x_i$ either to 1 or to 0 based on the structure of the starting graph $\mathcal{G}$.

The first type of inequalities, defined \textit{degree-cut inequalities}, explicitly impose that, for each node of the graph, at least one incident edge should be selected in the solution:
\begin{equation}\label{eq:degree-cut}
    \sum\limits_{i: e_i \in N(u)} x_i \geq 1 \quad \forall u \in V.
\end{equation}
When solving an LP relaxation of the original formulation, these constraints may help in finding a better bound of the optimal value. 
These $n$ inequalities are added to the MSTC model~\eqref{eq:mstc} as a-priori constraints.

The \textit{conflict-cycle inequalities} are a stronger version of the subtour elimination constraints~\eqref{eq:mstc2} obtained by taking into account also the conflicts among the edges.
Indeed, given a cycle $H\subset E$, constraints~\eqref{eq:mstc2} impose that the edges belonging to $P$ selected in any feasible solution of the problem should be at most $|H|-1$ (not all of them can be selected, otherwise a subtour is included in the solution). If an edge $e_c$ outside the cycle $H$ is in conflict with two edges belonging to the cycle $H,$ then in any feasible solution of MSTC, the number of selected edges in $H \cup \{e_c\}$ should be still less than $|H|-1$: no cycle, as well as no pairs of edges in conflict, should be included in the solution. This is imposed by the following inequalities: 
\begin{equation}\label{eq:conflictcycle-cut}
    \sum\limits_{i: e_i \in H} x_i + x_c \leq |H| - 1 \quad \forall \text{ cycle } H\subset E, e_c \in E\setminus H.
\end{equation}
These inequalities, as well as the classical subtour elimination constraints~\eqref{eq:mstc2}, are exponentially many, thus a separation procedure should be implemented for them. In particular, the separation problem for constraints~\eqref{eq:mstc2}, which consists in finding the violated constraints of type~\eqref{eq:mstc2} by a given solution, corresponds to solving a maximum-flow (minimum-cut) problem in a correspondingly defined capacitated graph as described in \cite{padberg}. As regards conflict-cycle inequalities~\eqref{eq:conflictcycle-cut}, the heuristic separation procedure proposed in \cite{carrabs2021bc} is implemented.

\paragraph{Preprocessing procedure}\label{sec:preprocess}
To the graph $\mathcal{G}$ taken into account, we apply the preprocessing procedure proposed in \cite{samer2015}. It consists of three steps, which are executed iteratively: each phase is executed as long as the problem instance is updated. Out of these two steps, we only perform the first two, as the last step was the most time-consuming and not particularly effective.

The first step looks for the so-called \textit{bridges} in $\mathcal{G}$, i.e., the edges the deletion of which disconnects the graph $\mathcal{G}$. Any bridge $e_i$, for its own nature, belongs to all the spanning trees of the graph. Thus, any feasible solution of MSTC cannot contain edges $e_j$ such that $(e_i,e_j) \in C$. These edges can therefore be removed from $\mathcal{G}$. 

The second step consists in selecting an edge $e_i$ and verifying if, due to the fact that edges $e_j$ such that $(e_i,e_j) \in C$ cannot be selected, the resulting graph is disconnected. If this is the case, then edge $e_i$ is removed from the graph. 
If any edge is removed during this phase, then the procedure restarts from the first one.

\section{Computational Tests}\label{sec:tests}
In this section, we report the results obtained by solving the instances proposed in \cite{zhang2011}\footnote{Instances available at \url{http://ic.ufal.br/professor/rian/ Instances/Instances-MCFST-Zhang.zip}.}, as well as the ones proposed in \cite{carrabs2021bc}\footnote{Instances available at \url{http://www.dipmat2.unisa.it/people/carrabs/www/DataSet/CMST.zip}.}, with \texttt{KS} (Algorithms~\ref{algo:kernel_search-if}--\ref{algo:kernel_search-ef}). In the following, we refer to 
the instances proposed in \cite{zhang2011} as \textit{ZKP instances} and to the ones proposed in \cite{carrabs2021bc} as \textit{CCPR instances}. 
The ZKP instances are generated by using different values for nodes, edges, and number of conflicts and they are classified into two types. By construction, there exists at least one conflict-free solution for all Type 2 instances, while this is not the case for Type 1 instances. In the computational tests, as for the Type 1, we consider only the 11 instances where at least the existence of a feasible solution is guaranteed.
The CCPR instances have a number of nodes in $\{25,50,75,100\}$ and a number of edges such that the density of the graph $m(m-1)/2$ is equal to $0.2$, $0.3$, or $0.4$. The number of conflicting pairs $|C|$ is equal to 1\%, 4\%, 7\% of the graph density. For each combination of these parameters (scenario), the authors generated 5 instances by varying the seed used for the generation of random numbers, obtaining a total of 180 instances. 

\subsection{Parameter tuning}\label{sec:tuning}
In this subsection, we explore the impact of various parameter settings on \texttt{KS} performance across different instances. 
We conducted an extensive parameter tuning process, exploring a range of values for the key parameters: $K$, $d$, $\Delta$, $P$. The goal is to identify the configurations that yield the best results in terms of both output quality and execution time.
To conduct a comprehensive analysis, we set up a series of experiments where we systematically varied key parameters of the algorithm. Each configuration is evaluated on multiple instances, allowing us to observe how different parameter settings influence the performance.

For the heuristic procedure, based on preliminary experiments, we set $\bar h = 20$ and $\bar t = 500$. For the whole algorithm, we set a time limit of $3600$ seconds and we test a time limit for each subproblem of $180, 300, 420$ seconds.

As for the dimension of the starting kernel set $K$ and the dimension of the starting buckets $d$, we need to take into account the dimension of the graph, in such a way that the starting kernel set contains a sufficient number of variables to eventually provide a feasible solution to the problem (i.e., $n-1$) plus possible alternatives, and the subproblems solved at each iteration over the kernel set and the buckets are not excessively large. For this reason, we considered the following formulas 
$$K = \lfloor \alpha\cdot(|V|-1) \rceil \quad d = \lfloor \beta\cdot(|E|-K)\rceil,$$ 
with $\alpha$ varying in $\{1.05,1.1,1.15,1.2,1.25\}$ and $\beta$ in $\{0.15,0.2,0.25\}$.
Instead, the maximum number of outer iterations $P$ is set to~$2, 3, 4$.
As for the parameter $\Delta$, indicating the maximum number of iterations \texttt{KS} can perform without improving, we need to link it to the number of buckets, that is, $b$. We set it to $\lfloor \delta \cdot b\rfloor$, with $\delta \in \{0.3,0.4,0.5,0.6\}$. We remark that $\Delta$ changes at every iteration $p$, since $b$ changes along these iterations.

In our preliminary experiments, we considered 9 Type 1 ZKP instances with up to 100 nodes and 500 edges and 10 CCPR instances with 50 nodes, 245 edges, and up to 1196 conflicts. We report, in Tables~\ref{tab:tuning}, for each parameter value the percentage of instances solved to optimality together with the average running time.
For the final tests, we select the configuration associated with the maximum weighted average between the percentage of instances for which we find a solution that is at least as good as the best found in the literature (with weight $0.8$) and the inverse of the average running time (with weight $0.2$). In the case of ties, we select the configuration with a greater percentage of optimally solved instances and then the one with a shorter time. For the ZKP instances, the selected parameter values are $P = 4, \alpha=1.1, \beta = 0.2, \delta = 0.6,$ and the inner time limit is $420$ seconds. For the CCPR instances, the selected parameter values are $P = 4, \alpha=1.2, \beta = 0.2, \delta =0.4,$ and the inner time limit is $180$ seconds.

\begin{table}[!ht]
    \centering
    \begin{minipage}{0.5\textwidth}
        \centering
        \subcaption{ZKP instances}
        \resizebox*{0.8\textwidth}{!}{  \begin{tabular}{|l|l|l|l|}
    \hline
        \rowcolor[HTML]{9B9B9B}parameter & \% solved instances & average time[s] & weighted average \\ \hline
        \rowcolor[HTML]{C0C0C0}$P$ & \multicolumn{3}{c|}{}  \\ \hline
        2 & 36.05 & 7.90 & 0.16\\ \hline
        3 & 47.84 & 14.90 & 0.20\\ \hline
        4 & 53.15 & 29.89 & 0.22\\ \hline
        \rowcolor[HTML]{C0C0C0}$\alpha$ & \multicolumn{3}{c|}{} \\ \hline
        1.05 & 40.12 & 19.69 & 0.17\\ \hline
        1.10 & 49.79 & 29.94 & 0.20\\ \hline
        1.15 & 49.59 & 20.55 & 0.20\\ \hline
        1.20 & 41.67 & 7.26 & 0.18\\ \hline
        1.25 & 47.22 & 14.99 & 0.20\\ \hline
        \rowcolor[HTML]{C0C0C0}$\beta$ & \multicolumn{3}{c|}{} \\ \hline
        0.15 & 46.11 & 19.89 & 0.19\\ \hline
        0.20 & 52.10 & 17.10 & 0.21\\ \hline
        0.25 & 38.83 & 20.05 & 0.16\\ \hline
        \rowcolor[HTML]{C0C0C0}$\delta$ & \multicolumn{3}{c|}{}\\ \hline
        0.30 & 45.43 & 20.51 & 0.19\\ \hline
        0.40 & 45.51 & 16.35 & 0.19\\ \hline
        0.50 & 45.84 & 17.90 & 0.19\\ \hline
        0.60 & 45.93 & 20.72 & 0.19\\ \hline
        \rowcolor[HTML]{C0C0C0}timelimit & \multicolumn{3}{c|}{} \\ \hline
        180 & 45.31 & 14.83 & 0.19\\ \hline
        300 & 45.62 & 18.21 & 0.19\\ \hline
        420 & 46.11 & 23.50 & 0.19\\ \hline
    \end{tabular}}
        \label{tab:tuning1}
    \end{minipage}%
    \hfill
    \begin{minipage}{0.5\textwidth}
        \centering
        \subcaption{CCPR instances}
        \resizebox*{0.8\textwidth}{!}{  \begin{tabular}{|l|l|l|l|}
    \hline
        \rowcolor[HTML]{9B9B9B}parameter & \% solved instances & average time[s] & weighted average \\ \hline
        \rowcolor[HTML]{C0C0C0}$P$ & \multicolumn{3}{c|}{}  \\ \hline
        2  & 61.17 & 2.26 & 0.29\\ \hline
        3  & 66.00 & 2.69 & 0.30\\ \hline
        4  & 66.00 & 2.68 & 0.30\\ \hline
        \rowcolor[HTML]{C0C0C0}$\alpha$ & \multicolumn{3}{c|}{} \\ \hline
        1.05  & 66.67 & 4.51 & 0.29\\ \hline
        1.10  & 63.33 & 3.57 & 0.28\\ \hline
        1.15  & 66.39 & 2.18 & 0.31\\ \hline
        1.20  & 65.56 & 0.81 & 0.39\\ \hline
        1.25  & 66.00 & 1.61 & 0.30\\ \hline
        \rowcolor[HTML]{C0C0C0}$\beta$ & \multicolumn{3}{c|}{} \\ \hline
        0.15 & 65.17 & 4.82 & 0.28\\ \hline
        0.20 & 64.00 & 1.31 & 0.33\\ \hline
        0.25 & 64.00 & 1.48 & 0.32\\ \hline
        \rowcolor[HTML]{C0C0C0}$\delta$ & \multicolumn{3}{c|}{}\\ \hline
        0.30 & 64.22 & 2.48 & 0.30\\ \hline
        0.40 & 64.44 & 2.56 & 0.30\\ \hline
        0.50 & 64.44 & 2.57 & 0.30\\ \hline
        0.60 & 64.44 & 2.58 & 0.30\\ \hline
        \rowcolor[HTML]{C0C0C0}timelimit & \multicolumn{3}{c|}{} \\ \hline
        180 & 64.39 & 2.50 & 0.30\\ \hline
        300 & 64.39 & 2.52 & 0.30\\ \hline
        420 & 64.39 & 2.51 & 0.30\\ \hline
    \end{tabular}}        \label{tab:tuning2}\end{minipage}
    \caption{Tuning of the parameters.}
    \label{tab:tuning}
\end{table}

\subsection{Comparisons with existing algorithms}
In this section, we report the results obtained by testing our \texttt{KS}, in comparison with state-of-the-art approaches.
Tables~\ref{tab:kernel_samur1}-\ref{tab:kernel_samur2} report the results on Type 1 and Type 2 ZKP instances (respectively) proposed in \cite{zhang2011} (for which at least a feasible solution is found), compared to the ones obtained by using the following heuristic approaches:
\begin{itemize}
    \item the Multi-Ethnic Genetic Approach (\texttt{MEGA}) proposed in \cite{carrabs2019};
    \item the Lagrangian-based heuristic \texttt{HDA+} proposed in \cite{gaudioso2021}, which provides a lower bound on the optimal value of the problem by solving the Lagrangian dual of the MSTC, and an upper bound through a local search procedure;
    \item the metaheuristic introduced in \cite{grasp}, \texttt{GRASP-AM}, which uses a GRASP enhanced with an adaptive memory strategy to efficiently explore the solution space and improve solution quality;
    \item the Two Phases Metaheuristic Technique (\texttt{TPMT}) proposed in \cite{chaubey2023}, which combines a hybrid artificial bee colony algorithm, and an iterated local search;
    \item the classic KS algorithm (\texttt{ClassicKS}) as introduced in \cite{Angelelli2010} -- never used to solve the MSTC problem -- which uses neither the independent set nor the heuristic. It performs at most $P$ iterations over the buckets for both types of instances, merging them two by two at each iteration. The stopping criterion is the same used in our KS algorithm, and we set $P=4, \alpha=1.2$, $\beta=0.1$, $\delta=0.3$, and the inner time limit to 420 seconds. These values are based on some preliminary experiments performed on a small subset of instances.
\end{itemize}
Following the approach used by \texttt{MEGA} \cite{carrabs2019}, \texttt{GRASP-AM} \cite{grasp}, \texttt{TMPT} \cite{chaubey2023} our proposed \texttt{KS} is executed for different runs on each instance (5 for the ZKP and 10 for the CCPR instances), each for a certain seed, which only affects the result found by the \texttt{StartingSolution} heuristic, used within the Initialization Phase.
Each subproblem involved in \texttt{KS} and \texttt{ClassicKS} algorithms is solved using the state-of-the-art solver Cplex.
They are coded in Python on an OSX platform, running on an Intel(R) Core(TM) i7-10700K CPU 3.80GHz with 32 GB of RAM, equipped with the IBM ILOG CPLEX (Version identifier: 22.1.0) solver. 
As concerns the time, in order to have a fair comparison, the CPU time reported in the literature for each method is properly scaled according to the Whetstone benchmark~\cite{whetstone2}.

The heading of Tables~\ref{tab:kernel_samur1}-\ref{tab:kernel_samur2} are the following: ID is the identifier of each instance; $n$ is the number of nodes of the graph; $m$ is the number of edges of the graph; $|C|$, is the number of conflicts; \textit{UB\_best} is the best value found either by the state-of-the-art methods we use to benchmark our approach or by the exact Branch and Cut method proposed in \cite{carrabs2021bc}; \textit{time}, the minimum computational time (in seconds) required by the method to return its best solution value; \textit{gap} which is the percentage gap of the best upper-bound $UB$ found over the multiple runs with respect to \textit{UB\_best}, i.e., \textit{gap} $=\frac{100(\text{\textit{UB}}-\text{\textit{UB\_best}})}{\text{\textit{UB\_best}}}.$
Whenever \textit{UB\_best} is certified to be optimal by the Branch and Cut, we report its value with the ``$*$'' symbol. Only in Table~\ref{tab:kernel_samur1}, for \texttt{KS}, we also report, in the \textit{UB\_time} column, the time needed by the algorithm to obtain the corresponding $UB$. We do not report this information in Table~\ref{tab:kernel_samur2} as \textit{UB\_time} corresponds to the whole time for the ZKP Type 2 instances. All the instances we solve are preprocessed according to the procedure described in \ref{sec:preprocess}. The time needed for this preprocessing as well as the time needed for all the algorithms involved within \texttt{KS} are included in the reported computational time. At the bottom of each table, the row ``Average values'' shows the average values for each column, and the row ``\%best'' reports the percentage of instances having a gap lower than or equal to 0\%. We highlight in bold the best average time and the best average gap.
\begin{table}[ht]
\centering
\resizebox*{0.95\textwidth}{!}{  
\begin{tabular}{|llll|c|rr|rr|rr|rr|rr|rrr|}
\rowcolor[HTML]{9B9B9B} 
\multicolumn{4}{|c|}{\cellcolor[HTML]{9B9B9B}Instance} & & \multicolumn{2}{c|}{\cellcolor[HTML]{9B9B9B}\texttt{MEGA}} & \multicolumn{2}{c|}{\cellcolor[HTML]{9B9B9B}\texttt{HDA+}} & \multicolumn{2}{c|}{\cellcolor[HTML]{9B9B9B}\texttt{GRASP-AM}} & \multicolumn{2}{c|}{\cellcolor[HTML]{9B9B9B}\texttt{TPMT}} &  \multicolumn{2}{c|}{\cellcolor[HTML]{9B9B9B}\texttt{ClassicKS}} & \multicolumn{3}{c|}{\cellcolor[HTML]{9B9B9B}\texttt{KS}}\\
\rowcolor[HTML]{C0C0C0} 
ID & $n$ & $m$ & $|C|$ & \textit{UB\_best} & \textit{time} & \textit{gap} & \textit{time} & \textit{gap} & \textit{time} & \textit{gap} & \textit{time} & \textit{gap} & \textit{time} & \textit{gap} & \textit{time} & \textit{UB\_time} & \textit{gap} \\
1 & 50 & 200 & 199 & 708$^*$ & 0.46 & 0.00 & 0.22 & 2.54 & 0.13 & 0.00 & 0.63 & 0.00 & 0.57 & 0.00 & 2.90 & 1.38 & 0.00 \\
2 & 50 & 200 & 398 & 770$^*$ & 0.44 & 0.00 & 0.29 & 1.30 & 0.13 & 0.00 & 0.79 & 0.00 & 16.42 & 0.00 & 1.24 & 0.57 & 0.00 \\
3 & 50 & 200 & 597 & 917$^*$ & 0.40 & 0.00 & 0.36 & 28.24 & 0.13 & 0.00 & 1.12 & 0.00 & 9.33 & 3.82 & 2.13 & 1.84 & 0.00 \\
4 & 50 & 200 & 995 & 1324$^*$ & 0.42 & 0.91 & 0.51 & 19.86 & 0.13 & 0.00 & 1.68 & 0.00 & 1.49 & 30.06 & 3.33 & 1.07 & 0.00 \\
5 & 100 & 300 & 448 & 4041$^*$ & 1.54 & 1.16 & 0.76 & 1.44 & 1.34 & 0.00 & 4.80 & 0.00 & 135.72 & 0.00 & 5.16 & 3.90 & 0.00 \\
6 & 100 & 300 & 897 & 5658$^*$ & 1.16 & 7.72 & 1.38 & 100 & 1.21 & 0.00 & 5.27 & 0.07 & 377.76 & 0.67 & 55.18 & 18.80 & 0.00 \\
7 & 100 & 500 & 1247 & 4275$^*$ & 3.33 & 0.00 & 1.90 & 0.61 & 1.30 & 0.00 & 4.68 & 0.00 & 207.63 & 0.16 & 10.70 & 3.32 & 0.00 \\
8 & 100 & 500 & 2495 & 5997$^*$ & 3.29 & 3.37 & 3.53 & 36.97 & 3.22 & 0.00 & 10.81 & 0.00 & 1016.28 & 1.62 & 75.45 & 26.35 & 0.00 \\
9 & 100 & 500 & 3741 & 7523 & 2.39 & 1.89 & 4.38 & 100 & 3.06 & 0.00 & 15.79 & 0.00 & 3985 & 100 & 308.82 & 24.77 & 0.00 \\
13 & 200 & 600 & 1797 & 13839 & 46.58 & 8.60 & 5.51 & 100 & 23.18 & 1.67 & 42.35 & 0.00 & 3881.71 & 100 & 1726.48 & 522.18 & -0.06 \\ 
16 & 200 & 800 & 3196 & 21518 & 67.62 & 2.75 & 12.17 & 100 & 38.68 & 0.89 & 29.36 & 0.00 & 1591.68 & 8.20 & 631.45 & 229.59 & -0.18 \\ 
\rowcolor[HTML]{C0C0C0} \multicolumn{5}{|c|}{Average values} & 11.60 & 2.40 & \textbf{2.82} & 44.63 & 6.59 & 0.23 & 10.66 & 0.01 & 1020.30 & 22.23 & 256.62 & 75.80 & \textbf{-0.02} \\ 
\rowcolor[HTML]{C0C0C0} \multicolumn{5}{|c|}{\%best} & \multicolumn{2}{c|}{36.4} & \multicolumn{2}{c|}{0} & \multicolumn{2}{c|}{81.8} & \multicolumn{2}{c|}{90.9} & \multicolumn{2}{c|}{27.3} & \multicolumn{3}{c|}{100} \\
    \end{tabular}}
\caption{Comparison of the results on feasible Type 1 ZKP instances proposed in \cite{zhang2011}.}
\label{tab:kernel_samur1}
\end{table}

As far as it concerns Type 1 instances (Table~\ref{tab:kernel_samur1}) the results of the average values line show that the fastest algorithm is \texttt{HDA+}, which requires only 2.82 seconds but is the less effective one with an average gap of 44.63\%. On the contrary, the most effective algorithm is \texttt{KS} which always finds the optimal solution, in the instances where it is known, and, for the last two instances, it reports a solution better than the best one known in the literature. Its final average gap is equal to -0.02\% and these results are obtained, on average, in around 75 seconds.

The effectiveness of \texttt{KS} on these instances is further certified by \%best line where its value is equal to 100\% since it always finds a solution that is at least as good as the best one known in the literature. Notice that it is the only algorithm to obtain this result. Similar results are obtained by \texttt{GRASP-AM} and \texttt{TPMT} which are faster than \texttt{KS} but show a higher average gap. Moreover, their \%best value is equal to 81.8\% and 90.91\%, respectively. 
It is worth noting the big difference between the results of \texttt{ClassicKS} and \texttt{KS}. \texttt{ClassicKS} is the slowest algorithm of the set and only in three cases it shows a gap equal to 0\%, while its average gap is very high (22.23\%). These preliminary results already show the impact that the use of the independent set strategy and the heuristic algorithm has on the effectiveness and performance of the KS method.

\begin{table}[ht]
\centering
\resizebox*{0.7\textwidth}{!}{  
\begin{tabular}{|llll|c|rr|rr|rr|rr|}
\rowcolor[HTML]{9B9B9B} 
\multicolumn{4}{|c|}{\cellcolor[HTML]{9B9B9B}Instance} & & \multicolumn{2}{c|}{\cellcolor[HTML]{9B9B9B}\texttt{MEGA}} & \multicolumn{2}{c|}{\cellcolor[HTML]{9B9B9B}\texttt{TPMT}} & \multicolumn{2}{c|}{\cellcolor[HTML]{9B9B9B}\texttt{ClassicKS}} & \multicolumn{2}{c|}{\cellcolor[HTML]{9B9B9B}\texttt{KS}}\\
\rowcolor[HTML]{C0C0C0} 
ID & $n$ & $m$ & $|C|$ & \textit{UB\_best} & \textit{time} & \textit{gap} & \textit{time} & \textit{gap} & \textit{time} & \textit{gap} & \textit{time} & \textit{gap} \\
24 & 50 & 200 & 3903 & 1636$^*$ & 0.30 & 0.00 & 0.04 & 0.00 & 0.02 & 0.00 & 0.02 & 0.00 \\
25 & 50 & 200 & 4877 & 2043$^*$ & 0.30 & 0.00 & 0.06 & 0.00 & 0.02 & 0.00 & 0.02 & 0.00 \\
26 & 50 & 200 & 5864 & 2338$^*$ & 0.33 & 0.00 & 0.06 & 0.00 & 0.02 & 0.00 & 0.02 & 0.00 \\
27 & 100 & 300 & 8609 & 7434$^*$ & 1.21 & 0.00 & 0.12 & 0.00 & 0.04 & 0.00 & 0.04 & 0.00 \\
28 & 100 & 300 & 10686 & 7968$^*$ & 1.08 & 0.00 & 0.16 & 0.00 & 0.02 & 0.00 & 0.02 & 0.00 \\
29 & 100 & 300 & 12761 & 8166$^*$ & 1.11 & 0.00 & 0.18 & 0.00 & 0.02 & 0.00 & 0.02 & 0.00 \\
30 & 100 & 500 & 24740 & 12652$^*$ & 2.12 & 0.00 & 0.63 & 0.00 & 0.31 & 0.00 & 0.31 & 0.00 \\
31 & 100 & 500 & 30886 & 11232$^*$ & 2.24 & 0.00 & 0.46 & 0.00 & 0.05 & 0.00 & 0.05 & 0.00 \\
32 & 100 & 500 & 36827 & 11481$^*$ & 2.26 & 0.00 & 0.45 & 0.00 & 0.06 & 0.00 & 0.06 & 0.00 \\
33 & 200 & 400 & 13660 & 17728$^*$ & 4.66 & 0.00 & 0.36 & 0.00 & 0.03 & 0.00 & 0.03 & 0.00 \\
34 & 200 & 400 & 17089 & 18617$^*$ & 4.81 & 0.00 & 0.37 & 0.00 & 0.03 & 0.00 & 0.03 & 0.00 \\
35 & 200 & 400 & 20470 & 19140$^*$ & 4.76 & 0.00 & 0.35 & 0.00 & 0.03 & 0.00 & 0.03 & 0.00 \\
36 & 200 & 600 & 34504 & 20716$^*$ & 7.18 & 0.00 & 1.05 & 0.00 & 0.08 & 0.00 & 0.08 & 0.00 \\
37 & 200 & 600 & 42860 & 18025$^*$ & 7.29 & 0.00 & 1.30 & 0.00 & 0.09 & 0.00 & 0.09 & 0.00 \\
38 & 200 & 600 & 50984 & 20864$^*$ & 7.96 & 0.00 & 1.24 & 0.00 & 0.10 & 0.00 & 0.10 & 0.00 \\
39 & 200 & 800 & 62625 & 39895$^*$ & 10.51 & 0.00 & 2.13 & 0.00 & 0.13 & 0.00 & 0.13 & 0.00 \\
40 & 200 & 800 & 78387 & 37671$^*$ & 10.09 & 0.00 & 1.91 & 0.00 & 0.14 & 0.00 & 0.14 & 0.00 \\
41 & 200 & 800 & 93978 & 38798$^*$ & 16.72 & 0.00 & 2.29 & 0.00 & 0.16 & 0.00 & 0.16 & 0.00 \\
42 & 300 & 600 & 31000 & 43721$^*$ & 11.96 & 0.00 & 0.80 & 0.00 & 0.05 & 0.00 & 0.05 & 0.00 \\
43 & 300 & 600 & 38216 & 44267$^*$ & 13.68 & 0.00 & 1.41 & 0.00 & 0.08 & 0.00 & 0.08 & 0.00 \\
44 & 300 & 600 & 45310 & 43071$^*$ & 15.63 & 0.00 & 1.28 & 0.00 & 0.07 & 0.00 & 0.07 & 0.00 \\
45 & 300 & 800 & 59600 & 43125$^*$ & 20.48 & 0.00 & 2.12 & 0.00 & 0.09 & 0.00 & 0.09 & 0.00 \\
46 & 300 & 800 & 74500 & 42292$^*$ & 22.05 & 0.00 & 2.40 & 0.00 & 0.11 & 0.00 & 0.11 & 0.00 \\
47 & 300 & 800 & 89300 & 44114$^*$ & 22.01 & 0.00 & 2.86 & 0.00 & 0.13 & 0.00 & 0.13 & 0.00 \\
48 & 300 & 1000 & 96590 & 71562$^*$ & 25.59 & 0.00 & 3.10 & 0.00 & 0.22 & 0.00 & 0.22 & 0.00 \\
49 & 300 & 1000 & 120500 & 76345$^*$ & 20.31 & 0.00 & 3.43 & 0.00 & 0.26 & 0.00 & 0.26 & 0.00 \\
50 & 300 & 1000 & 144090 & 78880$^*$ & 23.21 & 0.00 & 3.79 & 0.00 & 0.27 & 0.00 & 0.27 & 0.00 \\
\rowcolor[HTML]{C0C0C0} \multicolumn{5}{|c|}{Average values} & 9.62 & 0.00 & 1.27 & 0.00 & \textbf{0.10} & 0.00 & \textbf{0.10} & 0.00 \\
\rowcolor[HTML]{C0C0C0} \multicolumn{5}{|c|}{\%best} & \multicolumn{2}{c|}{100} & \multicolumn{2}{c|}{100} & \multicolumn{2}{c|}{100} & \multicolumn{2}{c|}{100} \\ 
\end{tabular}
}
\caption{Comparison of the results on ZKP Type 2 instances proposed in \cite{zhang2011}.}
    \label{tab:kernel_samur2}
\end{table}

As regards Type 2 instances (Table~\ref{tab:kernel_samur2}), in the literature they are solved by \texttt{MEGA}, and \texttt{TPMT} algorithms only. 
For this reason, Table~\ref{tab:kernel_samur2} reports the results of these two algorithms as well as \texttt{ClassicKS} and \texttt{KS}.
This set of instances appears very easy to solve for all the algorithms because all of them find the optimal known solution in less than 0.10 seconds, on average.

\begin{figure}[ht]
\begin{center}
  \includegraphics[scale=0.5]{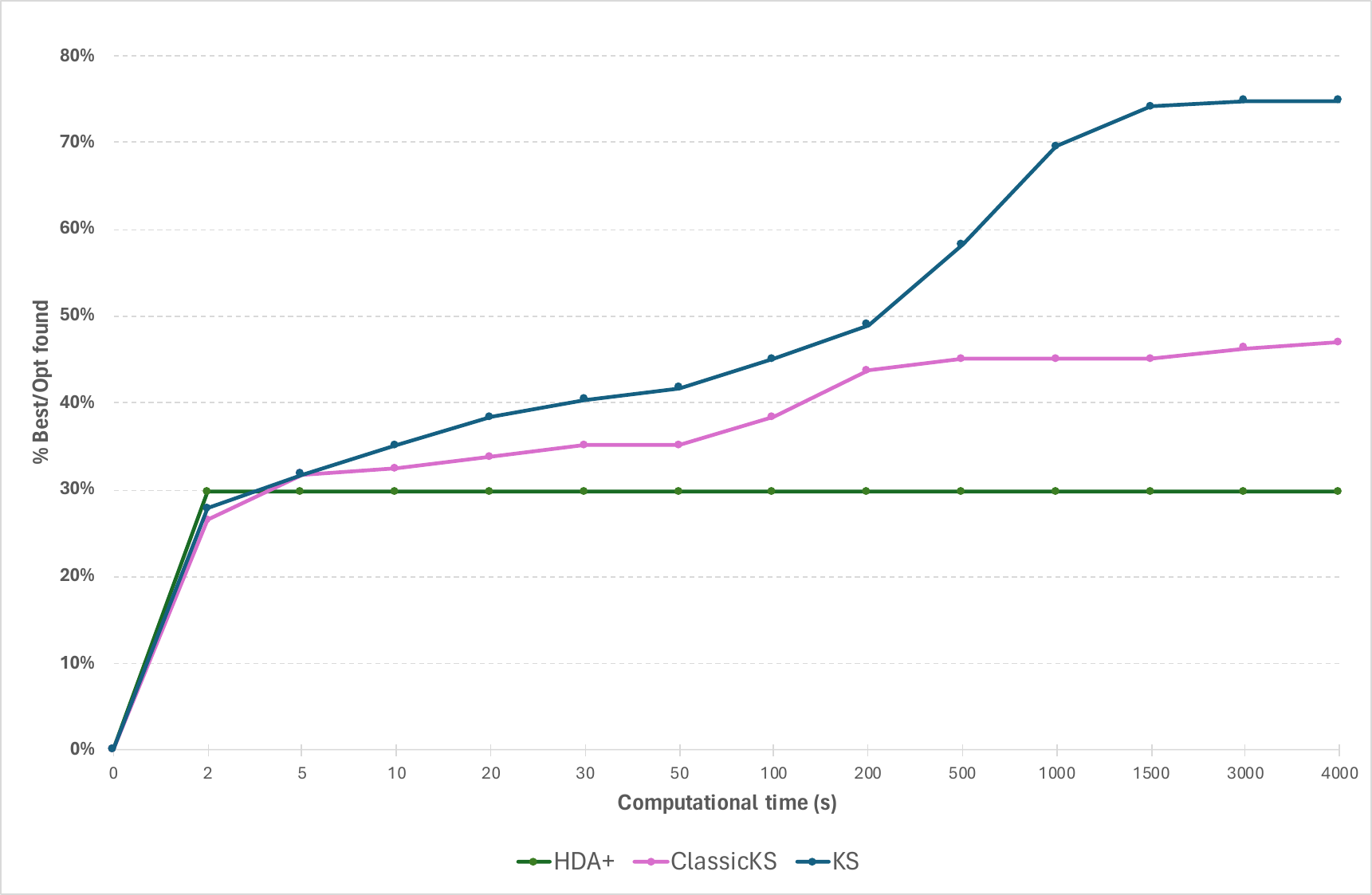}
\caption{Performance comparison among \texttt{HDA+}, \texttt{ClassicKS} and \texttt{KS} on the CCPR instances.}
\label{fig:chartL}
\end{center}
\end{figure}
\begin{figure}[ht]
\begin{center}
  \includegraphics[scale=0.42]{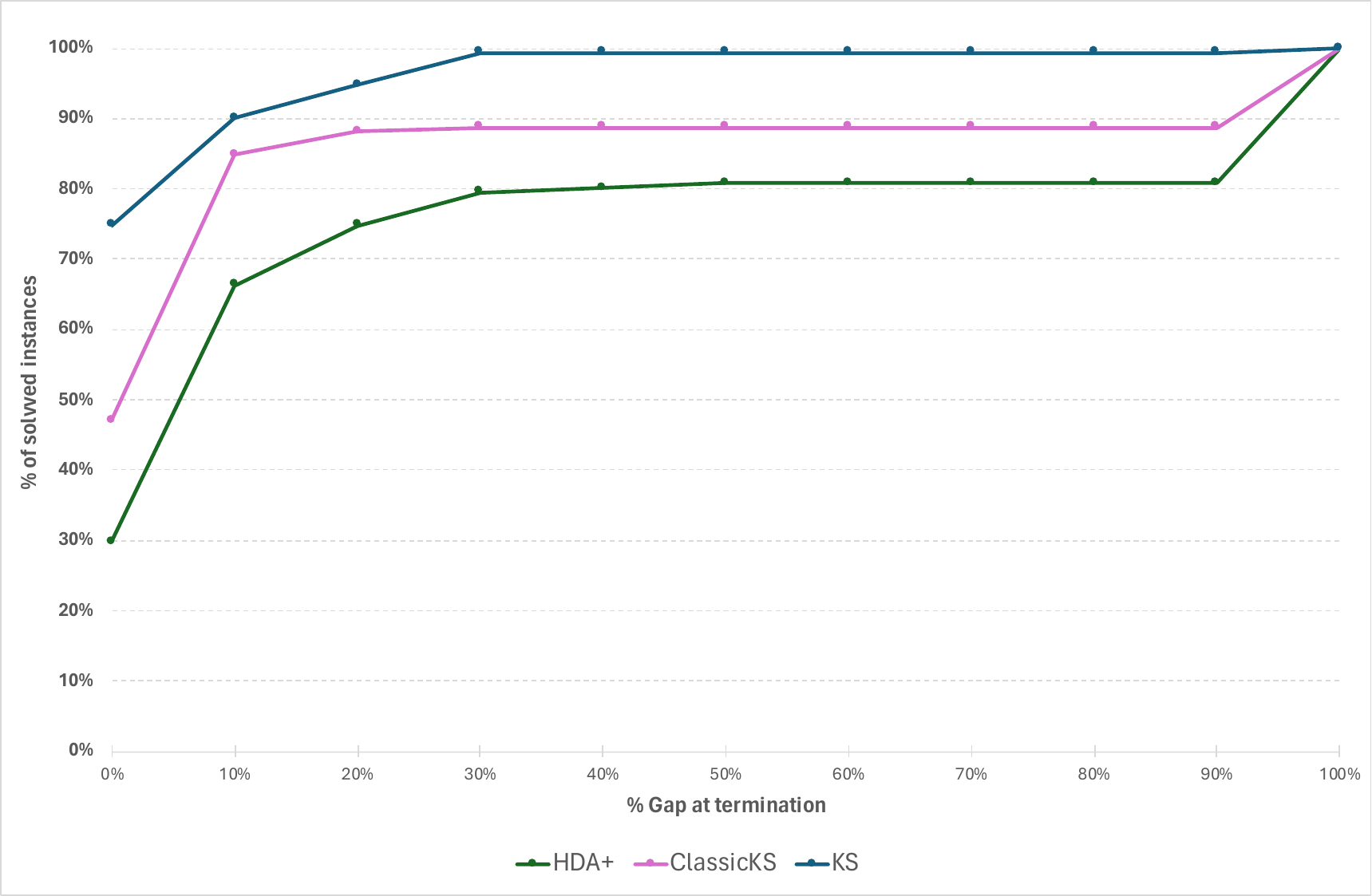}
\caption{Cumulative chart of the percentage of instances solved within a given gap at termination.}
\label{fig:chartM}
\end{center}
\end{figure}

The results on the CCPR instances are reported in Table~\ref{tab:kernel_cmst} in the Appendix. In this table, we only compare \texttt{KS} with \texttt{HDA+} method, besides \texttt{ClassicKS}, because the other algorithms were not tested on them.
As Table~\ref{tab:kernel_cmst} contains the results of 151 instances, we decided to move it into the Appendix and to represent its results both in Table~\ref{tab:kernel_cmst_aggregated} by aggregating them for the dimension of the instances, and in the cumulative charts reported in Figures~\ref{fig:chartL} and~\ref{fig:chartM}.

\begin{table}[ht]
\centering
\resizebox*{0.70\textwidth}{!}{  
\begin{tabular}{|ll|c|rr|rr|rrr|}
\rowcolor[HTML]{9B9B9B} 
\multicolumn{2}{|c|}{\cellcolor[HTML]{9B9B9B}Instance} & & \multicolumn{2}{c|}{\cellcolor[HTML]{9B9B9B}\texttt{HDA+}} & \multicolumn{2}{c|}{\cellcolor[HTML]{9B9B9B}\texttt{ClassicKS}} & \multicolumn{3}{c|}{\cellcolor[HTML]{9B9B9B}\texttt{KS}}\\
\rowcolor[HTML]{C0C0C0} 
$n$ & $m$ & \textit{UB\_best} & \textit{time} & \textit{gap} & \textit{time} & \textit{gap} & \textit{time} & \textit{UB\_time} & \textit{gap} \\
25 & 60 & 384.8 & 0.04 & 0.25 & 0.39 & 0.23 & 0.31 & 0.11 & 0.07 \\
25 & 90 & 329.9 & 0.07 & 0.29 & 1.60 & 0.43 & 0.77 & 0.26 & 0.00 \\
25 & 120 & 308.3 & 0.10 & 0.84 & 0.65 & 0.51 & 1.33 & 0.33 & 0.04 \\
50 & 245 & 704.6 & 0.61 & 34.9 & 616.69 & 8.12 & 73.79 & 11.8 & -0.39 \\
50 & 367 & 644.1 & 1.47 & 30.5 & 727.00 & 3.87 & 222.47 & 26.86 & -0.20 \\
50 & 490 & 600.8 & 2.87 & 26.54 & 800.76 & 3.48 & 309.19 & 42.17 & -0.32 \\
75 & 555 & 943.7 & 3.78 & 33.28 & 1189.01 & 5.28 & 609.41 & 87.99 & -0.04 \\
75 & 832 & 862.9 & 9.67 & 10.61 & 1572.92 & 0.83 & 707.01 & 107.55 & 0.07 \\
75 & 1110 & 862.4 & 21.92 & 19.68 & 1797.31 & 10.18 & 867.27 & 85.20 & 10.75\\
100 & 990 & 1285.7 & 17.41 & 51.46 & 2001.53 & 50.02 & 1084.52 & 113.25 & -3.37 \\
100 & 1485 & 1197.7 & 48.19 & 51.43 & 1988.07 & 50.13 & 1182.62 & 93.19 & 9.16 \\
100 & 1980 & 1124.4 & 102.92 & 39.53 & 2004.51 & 50.00 & 1207.13 & 95.62 & 9.71 \\
\rowcolor[HTML]{C0C0C0} \multicolumn{3}{|c|}{Average values} &  \textbf{14.16} & 23.05 & 924.09 & 12.72  & 440.83 & 47.26 & \textbf{1.73} \\ 
\rowcolor[HTML]{C0C0C0} \multicolumn{3}{|c|}{\%best} & \multicolumn{2}{c|}{29.80} & \multicolumn{2}{c|}{47.02} & \multicolumn{3}{c|}{74.83}\\
\hline
    \end{tabular}
}
\caption{Comparison of the aggregated results on CCPR instances proposed in \cite{carrabs2021bc}.}
\label{tab:kernel_cmst_aggregated}
\end{table}

Already by looking at the aggregated results in Table~\ref{tab:kernel_cmst_aggregated} it is evident that, despite requiring more computational time w.r.t.\ \texttt{HDA+}, \texttt{KS} is by far the best in terms of the gap. Specifically, \texttt{KS} achieves the best average gap across all instance groups and the highest percentage of instances with a non-positive gap relative to the best solution value found in the literature (last row of the table). Notably, 74.83\% of the CCPR instances are solved with a non-positive gap (15 new best-known solution), the average gap is 1.73\% and the average computing time is approximately 7 minutes. Although this is significantly higher than the average computing time of the \texttt{HDA+} method, the latter yields a much larger average gap of 23.05\%. The classic KS approach is slower than our KS, returning an average gap of 12.72\%, which is still less than the one returned by the \texttt{HDA+} method.

In Figure~\ref{fig:chartL}, the horizontal axis reports the computational time in seconds, and the vertical one shows the percentage of best/optimal solutions found within that time. More precisely, a $(x,y)$ point on this plot means that within $x$ seconds of computation, the algorithm found the best/optimal solution for y\% of the instances in at most $x$ seconds. This implies that the faster the growth of a curve, the better the performance. It is worth noting that in our tables we compare the solutions of our algorithms with the best ones available in the literature (i.e. \textit{UB\_best} value). As a consequence whenever \texttt{KS} finds a solution which is better than \textit{UB\_best}, the corresponding percentage gap will be negative. 
Since negative percentages are not allowed in the graphic of Figure~\ref{fig:chartL} then the comparison is not carried out with \textit{UB\_Best} but with the best solution obtained by all the algorithms of the literature as well as our algorithms. We refer to this new value as \textit{New\_UB\_Best}, and we will use it also for the construction of the graphic in Figure~\ref{fig:chartM}.
The blue curve in Figure~\ref{fig:chartL} is associated with \texttt{KS}, the pink one corresponds to \texttt{ClassicKS}, and the green one to the \texttt{HDA+} algorithm. Figure~\ref{fig:chartL} certifies the superiority of \texttt{KS} with respect to the other two algorithms. Indeed, \texttt{KS} reaches the $\sim$75\% of best/optimal solutions found, while for \texttt{ClassicKS} and \texttt{HDA+} this percentage is equal to $\sim$47\% and $\sim$30\% only, respectively. By comparing the \texttt{ClassicKS} and \texttt{KS} curves, we observe that \texttt{ClassicKS} requires the whole time limit to reach its peak of $\sim$47\% whereas \texttt{KS} reaches this percentage in less than 200 seconds. 
\texttt{HDA+} in total only manages to find 46 best solutions (see Figure~\ref{fig:chartM}) and it does it in 2 seconds (then it always obtains a positive gap at termination, that is why the curve stays constant). \texttt{KS} obtains a very close percentage (the 28\%) in the same amount time.

Figure~\ref{fig:chartM} shows the effectiveness of the three algorithms in the 151 instances considered. The horizontal axis reports the percentage gap from \textit{New\_UB\_Best} at the end of the computation, whereas the vertical one shows the percentage of instances for which the percentage gap returned by the algorithm at termination is lower than or equal to the value reported on the horizontal axis. As for the graphic in Figure~\ref{fig:chartL}, the faster the growth of a curve, the better the effectiveness of the corresponding algorithm. 
It is easy to see from this graphic that \texttt{KS} overcomes the other two algorithms because its curve is always higher than the other two ones. This means that \texttt{KS} obtains a lower percentage gap from the best solution, for more instances with respect to the other two algorithms.
When the horizontal value is 0\%, the curves of Figure~\ref{fig:chartM} show the percentage of instances for which the solution of the algorithm matches the \textit{New\_UB\_Best} value, and this percentage is equal to $\sim$75\%, $\sim$47\% and $\sim$30\% for \texttt{KS}, \texttt{ClassicKS} and \texttt{HDA+}, respectively. These percentages clearly coincide with the peaks reported in Figure~\ref{fig:chartL} for these three algorithms. The percentage of instances for which the gap is lower than 10\% significantly increases for all the algorithms. This percentage is, indeed, equal to $\sim$90\%, $\sim$85\% and $\sim$66\% for \texttt{KS}, \texttt{ClassicKS} and \texttt{HDA+}, respectively. These results highlight the effectiveness of \texttt{KS} which almost always returns a solution very close to the best one and obtains a gap lower than or equal to 20\% for the 95\% of instances. Within the same threshold of 20\%, the percentage of \texttt{ClassicKS} and \texttt{HDA+} equals $\sim$88\% and $\sim$75\%, respectively, which are much lower than the one of \texttt{KS}.

In conclusion, while the average termination time of our KS algorithm is always higher than that of other algorithms (except for the classic KS one, where our enhancements lead to improvements also in terms of time), the time required to find the best-returned incumbent is definitely competitive. Moreover, this best solution is 17 times better than those reported in the literature and at least as good as the best-known solutions in 80\% of cases.

\section{Conclusion}\label{sec:conclusion}
In this paper, we introduced a novel solution approach to the Minimum Spanning Tree with Conflicts problem, leveraging the KS method \cite{Angelelli2010}. This algorithm focuses on obtaining a feasible solution, ideally of high quality, by strategically selecting a small set of promising variables. This subset, known as the kernel set, is initially constructed through the solution of the LP relaxation of the original ILP. Throughout the algorithm, the kernel set is dynamically enriched by identifying new promising variables through the solution of ILP subproblems of increasing size.
A key innovation in our methodology lies in the incorporation of independent sets of proper subgraphs within the KS algorithm. This novel approach draws inspiration from the observation that the edges of any feasible MSTC solution correspond to an independent set in the conflict graph. The independent sets are used in order to both detect affinities among variables and find a starting point for the algorithm.

Through comparative analysis with existing heuristics, we demonstrated the efficacy of our method in producing competitive results. It is evident that our tailored KS approach offers a promising tool for tackling the MSTC problem and potentially other conflict-based optimization challenges. 

\bibliographystyle{plain}
\bibliography{MSTC-kernel}

\appendix
\setcounter{table}{0}
\renewcommand{\thetable}{A\arabic{table}}
\section{CCPR results}\label{appendix}
In this section, we show the detailed results of the comparison among the \texttt{HDA+}, \texttt{ClassicKS}, and \texttt{KS} on the CCPR instances proposed in \cite{carrabs2021bc}. As already noted for the other tables, whenever the B\&C method finds the optimal value, we report such value with the ``*'' symbol in the \textit{UB\_best} column. Whenever no feasible solution is known in the literature, but we do find one, a ``-'' is reported in the \textit{UB\_best} column. In the \textit{gap} columns, if the considered method does not find a feasible solution, we report a $100$. Instead, if the considered method does find a feasible solution, but there is a ``-'' in the \textit{UB\_best} column, we report a $-100$. Also in this case like for the previous tables, we report the results only on the instances where at least one of the methods (including ours) finds a feasible solution. As a consequence, the instances with ID 151--155, 166--170, 181, 182, 184, 185, 196--200, 211--215, 226--230 are not present in this table.
The headings of Table~\ref{tab:kernel_cmst} are the same as in Table~\ref{tab:kernel_samur1}. Again, we report the ``Average values'' and the ``\%best'' rows at the end of the table. 

{\fontsize{6}{10} \selectfont\begin{longtable}{!{\color{black!20}\vrule width 0.8pt}
llll!{\color{black!20}\vrule width 0.8pt}c!{\color{black!20}\vrule width  0.8pt}rr!{\color{black!20}\vrule width  0.8pt}rr!{\color{black!20}\vrule width  0.8pt}rrr!{\color{black!20}\vrule width  0.8pt}}
\caption{Comparison of the results on the CCPR instances proposed in \cite{carrabs2021bc}, for which at least a feasible solution is found.}
\endlastfoot
\rowcolor[HTML]{9B9B9B} 
\multicolumn{4}{!{\color{black!20}\vrule width  0.8pt}c!{\color{black!20}\vrule width  0.8pt}}{\cellcolor[HTML]{9B9B9B}Instance} &  & \multicolumn{2}{c!{\color{black!20}\vrule width  0.8pt}}{\cellcolor[HTML]{9B9B9B}\texttt{HDA+}} & \multicolumn{2}{c!{\color{black!20}\vrule width  0.8pt}}{\cellcolor[HTML]{9B9B9B}{ClassicKS}} & \multicolumn{3}{c!{\color{black!20}\vrule width  0.8pt}}{\cellcolor[HTML]{9B9B9B}KS} \\
\rowcolor[HTML]{C0C0C0} 
ID & $n$ & $m$ & $|C|$ & \textit{UB\_best} & \textit{time} & \textit{gap} & \textit{time} & \textit{gap} & \textit{time} & \textit{UB\_time} & \textit{gap} \\
\endhead
\rowcolor{black!10} 51 & 25 & 60 & 18 & 347$^*$ & 0.03 & 0.00 & 0.02 & 0.00 & 0.02 & 0.02 & 0.00 \\ 
 52 & 25 & 60 & 18 & 389$^*$ & 0.03 & 0.00 & 0.00 & 0.00 & 0.00 & 0.00 & 0.00 \\ 
\rowcolor{black!10} 53 & 25 & 60 & 18 & 353$^*$ & 0.03 & 0.00 & 4.15 & 1.98 & 3.16 & 0.95 & 0.00 \\ 
 54 & 25 & 60 & 18 & 346$^*$ & 0.03 & 0.00 & 0.00 & 0.00 & 0.00 & 0.00 & 0.00 \\ 
\rowcolor{black!10} 55 & 25 & 60 & 18 & 336$^*$ & 0.03 & 0.00 & 0.00 & 0.00 & 0.00 & 0.00 & 0.00 \\ 
 56 & 25 & 60 & 71 & 381$^*$ & 0.04 & 0.00 & 0.42 & 0.00 & 0.15 & 0.11 & 0.00 \\ 
\rowcolor{black!10} 57 & 25 & 60 & 71 & 390$^*$ & 0.04 & 0.00 & 0.62 & 0.26 & 0.77 & 0.09 & 0.26 \\ 
 58 & 25 & 60 & 71 & 372$^*$ & 0.04 & 0.00 & 0.02 & 0.00 & 0.04 & 0.04 & 0.00 \\ 
\rowcolor{black!10} 59 & 25 & 60 & 71 & 357$^*$ & 0.04 & 0.00 & 0.02 & 0.00 & 0.04 & 0.04 & 0.00 \\ 
 60 & 25 & 60 & 71 & 406$^*$ & 0.04 & 0.00 & 0.00 & 0.00 & 0.00 & 0.00 & 0.00 \\ 
\rowcolor{black!10} 61 & 25 & 60 & 124 & 385$^*$ & 0.04 & 0.00 & 0.00 & 0.00 & 0.00 & 0.00 & 0.00 \\ 
 62 & 25 & 60 & 124 & 432$^*$ & 0.04 & 0.00 & 0.00 & 0.00 & 0.00 & 0.00 & 0.00 \\ 
\rowcolor{black!10} 63 & 25 & 60 & 124 & 458$^*$ & 0.06 & 3.49 & 0.20 & 0.00 & 0.32 & 0.21 & 0.00 \\ 
 64 & 25 & 60 & 124 & 400$^*$ & 0.05 & 0.00 & 0.26 & 1.25 & 0.07 & 0.04 & 0.75 \\ 
\rowcolor{black!10} 65 & 25 & 60 & 124 & 420$^*$ & 0.08 & 0.24 & 0.12 & 0.00 & 0.12 & 0.08 & 0.00 \\ 
 66 & 25 & 90 & 41 & 311$^*$ & 0.03 & 0.00 & 2.80 & 0.00 & 0.25 & 0.10 & 0.00 \\ 
\rowcolor{black!10} 67 & 25 & 90 & 41 & 306$^*$ & 0.04 & 0.00 & 0.00 & 0.00 & 0.00 & 0.00 & 0.00 \\ 
 68 & 25 & 90 & 41 & 299$^*$ & 0.03 & 0.00 & 0.69 & 0.00 & 0.24 & 0.14 & 0.00 \\ 
\rowcolor{black!10} 69 & 25 & 90 & 41 & 297$^*$ & 0.03 & 0.00 & 5.08 & 4.38 & 0.75 & 0.47 & 0.00 \\ 
 70 & 25 & 90 & 41 & 318$^*$ & 0.03 & 0.00 & 12.87 & 0.00 & 5.95 & 0.79 & 0.00 \\ 
\rowcolor{black!10} 71 & 25 & 90 & 161 & 305$^*$ & 0.06 & 0.00 & 0.05 & 0.00 & 0.05 & 0.04 & 0.00 \\ 
 72 & 25 & 90 & 161 & 339$^*$ & 0.07 & 0.00 & 0.00 & 0.00 & 0.00 & 0.00 & 0.00 \\ 
\rowcolor{black!10} 73 & 25 & 90 & 161 & 344$^*$ & 0.06 & 0.00 & 0.02 & 0.00 & 0.04 & 0.06 & 0.00 \\ 
 74 & 25 & 90 & 161 & 329$^*$ & 0.07 & 0.61 & 0.10 & 0.30 & 0.24 & 0.20 & 0.00 \\ 
\rowcolor{black!10} 75 & 25 & 90 & 161 & 326$^*$ & 0.07 & 0.31 & 0.59 & 1.53 & 2.74 & 1.26 & 0.00 \\ 
 76 & 25 & 90 & 281 & 349$^*$ & 0.08 & 0.00 & 0.25 & 0.00 & 0.41 & 0.34 & 0.00 \\ 
\rowcolor{black!10} 77 & 25 & 90 & 281 & 385$^*$ & 0.10 & 0.00 & 0.91 & 0.00 & 0.22 & 0.14 & 0.00 \\ 
 78 & 25 & 90 & 281 & 335$^*$ & 0.09 & 0.00 & 0.17 & 0.30 & 0.20 & 0.12 & 0.00 \\ 
\rowcolor{black!10} 79 & 25 & 90 & 281 & 348$^*$ & 0.12 & 2.87 & 0.22 & 0.00 & 0.21 & 0.11 & 0.00 \\ 
 80 & 25 & 90 & 281 & 357$^*$ & 0.11 & 0.56 & 0.16 & 0.00 & 0.22 & 0.15 & 0.00 \\ 
\rowcolor{black!10} 81 & 25 & 120 & 72 & 282$^*$ & 0.04 & 0.00 & 0.00 & 0.00 & 0.00 & 0.00 & 0.00 \\ 
 82 & 25 & 120 & 72 & 294$^*$ & 0.04 & 0.00 & 0.97 & 0.00 & 15.40 & 3.44 & 0.00 \\ 
\rowcolor{black!10} 83 & 25 & 120 & 72 & 284$^*$ & 0.05 & 0.00 & 0.00 & 0.00 & 0.00 & 0.00 & 0.00 \\ 
 84 & 25 & 120 & 72 & 281$^*$ & 0.04 & 0.00 & 2.56 & 0.00 & 1.90 & 0.10 & 0.00 \\ 
\rowcolor{black!10} 85 & 25 & 120 & 72 & 292$^*$ & 0.04 & 0.00 & 0.00 & 0.00 & 0.00 & 0.00 & 0.00 \\ 
 86 & 25 & 120 & 286 & 321$^*$ & 0.10 & 0.00 & 0.14 & 0.31 & 0.23 & 0.11 & 0.00 \\ 
\rowcolor{black!10} 87 & 25 & 120 & 286 & 317$^*$ & 0.11 & 0.00 & 0.51 & 0.00 & 0.21 & 0.08 & 0.00 \\ 
 88 & 25 & 120 & 286 & 284$^*$ & 0.08 & 0.00 & 0.00 & 0.00 & 0.00 & 0.00 & 0.00 \\ 
\rowcolor{black!10} 89 & 25 & 120 & 286 & 311$^*$ & 0.11 & 0.32 & 0.07 & 0.00 & 0.03 & 0.03 & 0.00 \\ 
 90 & 25 & 120 & 286 & 290$^*$ & 0.08 & 0.00 & 0.02 & 0.00 & 0.06 & 0.06 & 0.00 \\ 
\rowcolor{black!10} 91 & 25 & 120 & 500 & 329$^*$ & 0.17 & 3.65 & 0.71 & 1.82 & 0.24 & 0.19 & 0.00 \\ 
 92 & 25 & 120 & 500 & 339$^*$ & 0.16 & 2.36 & 2.83 & 0.00 & 0.63 & 0.48 & 0.00 \\ 
\rowcolor{black!10} 93 & 25 & 120 & 500 & 368$^*$ & 0.19 & 4.08 & 1.50 & 0.82 & 0.31 & 0.18 & 0.00 \\ 
 94 & 25 & 120 & 500 & 311$^*$ & 0.14 & 0.97 & 0.14 & 0.00 & 0.38 & 0.11 & 0.00 \\ 
\rowcolor{black!10} 95 & 25 & 120 & 500 & 321$^*$ & 0.16 & 1.25 & 0.19 & 0.62 & 0.57 & 0.12 & 0.62 \\ 
 96 & 50 & 245 & 299 & 619$^*$ & 0.24 & 0.00 & 4.27 & 0.00 & 442.27 & 71.54 & 0.00 \\ 
\rowcolor{black!10} 97 & 50 & 245 & 299 & 604$^*$ & 0.24 & 0.00 & 0.29 & 0.00 & 57.32 & 13.16 & 0.00 \\ 
 98 & 50 & 245 & 299 & 634$^*$ & 0.22 & 0.00 & 0.01 & 0.00 & 0.01 & 0.01 & 0.00 \\ 
\rowcolor{black!10} 99 & 50 & 245 & 299 & 616$^*$ & 0.26 & 0.00 & 0.35 & 0.16 & 6.72 & 2.64 & 0.00 \\ 
 100 & 50 & 245 & 299 & 595$^*$ & 0.28 & 0.00 & 21.59 & 0.00 & 0.25 & 0.25 & 0.00 \\ 
\rowcolor{black!10} 101 & 50 & 245 & 1196 & 678$^*$ & 0.62 & 2.95 & 10.90 & 0.00 & 24.81 & 5.09 & 0.00 \\ 
 102 & 50 & 245 & 1196 & 681$^*$ & 0.64 & 5.87 & 86.60 & 3.08 & 26.29 & 8.49 & 0.00 \\ 
\rowcolor{black!10} 103 & 50 & 245 & 1196 & 709$^*$ & 0.64 & 2.26 & 17.80 & 0.71 & 12.86 & 3.83 & 0.00 \\ 
 104 & 50 & 245 & 1196 & 639$^*$ & 0.55 & 2.66 & 11.05 & 0.16 & 6.51 & 2.19 & 0.00 \\ 
\rowcolor{black!10} 105 & 50 & 245 & 1196 & 681$^*$ & 0.65 & 9.84 & 23.49 & 0.59 & 15.44 & 4.55 & 0.29 \\ 
 106 & 50 & 245 & 2093 & 833 & 1.00 & 100.00 & 1102.05 & 1.20 & 149.16 & 9.95 & -2.04 \\ 
\rowcolor{black!10} 107 & 50 & 245 & 2093 & 835 & 1.02 & 100.00 & 3685.50 & 100.00 & 201.32 & 10.71 & 0.00 \\ 
 108 & 50 & 245 & 2093 & 840 & 0.91 & 100.00 & 2245.91 & 7.50 & 79.02 & 9.81 & -3.57 \\ 
\rowcolor{black!10} 109 & 50 & 245 & 2093 & 836 & 0.97 & 100.00 & 1721.24 & 5.50 & 34.16 & 17.88 & -0.48 \\ 
 110 & 50 & 245 & 2093 & 769$^*$ & 0.89 & 100.00 & 319.20 & 2.86 & 50.74 & 16.93 & 0.00 \\ 
\rowcolor{black!10} 111 & 50 & 367 & 672 & 570$^*$ & 0.57 & 0.00 & 0.49 & 0.35 & 105.67 & 21.60 & 0.00 \\ 
 112 & 50 & 367 & 672 & 561$^*$ & 0.63 & 0.00 & 4.10 & 0.00 & 750.19 & 77.86 & 0.00 \\ 
\rowcolor{black!10} 113 & 50 & 367 & 672 & 573$^*$ & 0.68 & 0.00 & 0.53 & 0.00 & 8.60 & 2.73 & 0.00 \\ 
 114 & 50 & 367 & 672 & 560$^*$ & 0.56 & 0.00 & 0.09 & 0.00 & 0.88 & 0.66 & 0.00 \\ 
\rowcolor{black!10} 115 & 50 & 367 & 672 & 549$^*$ & 0.64 & 0.36 & 0.44 & 0.18 & 203.51 & 0.29 & 0.00 \\ 
 116 & 50 & 367 & 2687 & 612$^*$ & 1.52 & 7.35 & 34.96 & 0.33 & 277.78 & 44.35 & 0.00 \\ 
\rowcolor{black!10} 117 & 50 & 367 & 2687 & 615$^*$ & 1.44 & 7.81 & 34.21 & 0.81 & 19.76 & 4.04 & 0.00 \\ 
 118 & 50 & 367 & 2687 & 587$^*$ & 1.50 & 8.18 & 23.60 & 1.53 & 2.46 & 0.97 & 0.00 \\ 
\rowcolor{black!10} 119 & 50 & 367 & 2687 & 634$^*$ & 1.54 & 13.72 & 34.00 & 1.74 & 86.84 & 16.41 & 0.00 \\ 
 120 & 50 & 367 & 2687 & 643$^*$ & 1.57 & 7.00 & 15.72 & 2.33 & 8.08 & 3.19 & 0.00 \\ 
\rowcolor{black!10} 121 & 50 & 367 & 4702 & 726 & 2.26 & 100.00 & 3400.73 & 13.36 & 632.18 & 72.62 & 1.10 \\ 
 122 & 50 & 367 & 4702 & 770 & 2.19 & 100.00 & 2418.38 & 16.36 & 116.46 & 15.51 & -1.30 \\ 
\rowcolor{black!10} 123 & 50 & 367 & 4702 & 786 & 2.30 & 100.00 & 2982.65 & 7.25 & 483.81 & 57.09 & -2.54 \\ 
 124 & 50 & 367 & 4702 & 711 & 2.42 & 100.00 & 952.21 & 4.92 & 280.06 & 58.20 & -1.97 \\ 
\rowcolor{black!10} 125 & 50 & 367 & 4702 & 764 & 2.16 & 13.61 & 1002.78 & 8.90 & 360.74 & 27.39 & 1.70 \\ 
 126 & 50 & 490 & 1199 & 548$^*$ & 1.30 & 0.73 & 29.42 & 0.00 & 118.80 & 14.62 & 0.00 \\ 
\rowcolor{black!10} 127 & 50 & 490 & 1199 & 530$^*$ & 1.21 & 0.19 & 214.98 & 0.00 & 883.20 & 50.26 & 0.00 \\ 
 128 & 50 & 490 & 1199 & 549$^*$ & 1.13 & 0.00 & 0.12 & 0.00 & 181.96 & 19.68 & 0.00 \\ 
\rowcolor{black!10} 129 & 50 & 490 & 1199 & 540$^*$ & 1.29 & 0.19 & 0.10 & 0.00 & 0.09 & 0.08 & 0.00 \\ 
 130 & 50 & 490 & 1199 & 540$^*$ & 1.10 & 0.00 & 56.15 & 0.00 & 120.56 & 0.51 & 0.19 \\ 
\rowcolor{black!10} 131 & 50 & 490 & 4793 & 594$^*$ & 2.98 & 5.89 & 40.07 & 0.34 & 131.50 & 26.89 & 0.00 \\ 
 132 & 50 & 490 & 4793 & 579$^*$ & 3.15 & 12.26 & 69.87 & 0.00 & 23.29 & 3.00 & 0.00 \\ 
\rowcolor{black!10} 133 & 50 & 490 & 4793 & 589$^*$ & 3.05 & 11.55 & 17.78 & 0.17 & 4.69 & 2.80 & 0.00 \\ 
 134 & 50 & 490 & 4793 & 577$^*$ & 3.05 & 11.44 & 29.16 & 2.25 & 4.88 & 2.87 & 0.00 \\ 
\rowcolor{black!10} 135 & 50 & 490 & 4793 & 592$^*$ & 3.03 & 13.18 & 36.84 & 0.51 & 71.11 & 13.86 & 0.00 \\ 
 136 & 50 & 490 & 8387 & 678 & 4.20 & 19.76 & 2445.19 & 16.96 & 476.13 & 48.62 & -1.18 \\ 
\rowcolor{black!10} 137 & 50 & 490 & 8387 & 651 & 4.13 & 100.00 & 3606.37 & 15.36 & 998.68 & 155.43 & -0.61 \\ 
 138 & 50 & 490 & 8387 & 689 & 4.27 & 100.00 & 1893.81 & 1.16 & 396.10 & 76.11 & -1.45 \\ 
\rowcolor{black!10} 139 & 50 & 490 & 8387 & 682 & 4.57 & 100.00 & 1301.15 & 9.09 & 290.68 & 42.49 & 0.73 \\ 
 140 & 50 & 490 & 8387 & 674 & 4.47 & 22.85 & 2270.21 & 6.38 & 936.19 & 175.34 & -2.52 \\ 
\rowcolor{black!10} 141 & 75 & 555 & 1538 & 868$^*$ & 1.96 & 0.12 & 25.20 & 0.23 & 573.26 & 30.45 & 0.00 \\ 
 142 & 75 & 555 & 1538 & 871$^*$ & 2.28 & 0.80 & 3.38 & 0.00 & 772.83 & 99.49 & 0.00 \\ 
\rowcolor{black!10} 143 & 75 & 555 & 1538 & 838$^*$ & 2.33 & 0.72 & 0.84 & 0.00 & 19.99 & 7.48 & 0.00 \\ 
 144 & 75 & 555 & 1538 & 855$^*$ & 1.79 & 0.00 & 18.13 & 0.12 & 638.29 & 104.89 & 0.00 \\ 
\rowcolor{black!10} 145 & 75 & 555 & 1538 & 857$^*$ & 2.01 & 0.23 & 0.96 & 0.23 & 10.30 & 1.32 & 0.00 \\ 
 146 & 75 & 555 & 6150 & 1047 & 5.64 & 18.05 & 2511.64 & 6.21 & 400.59 & 43.88 & 0.48 \\ 
\rowcolor{black!10} 147 & 75 & 555 & 6150 & 1069 & 5.57 & 12.91 & 2310.07 & 6.08 & 1033.18 & 216.67 & -2.15 \\ 
 148 & 75 & 555 & 6150 & 1040 & 5.50 & 100.00 & 3805.55 & 20.87 & 754.06 & 95.82 & -0.58 \\ 
\rowcolor{black!10} 149 & 75 & 555 & 6150 & 998 & 5.05 & 100.00 & 1086.86 & 1.10 & 852.15 & 128.74 & 1.40 \\ 
 150 & 75 & 555 & 6150 & 994 & 5.69 & 100.00 & 2127.32 & 18.01 & 1039.43 & 151.19 & 0.40 \\ 
\rowcolor{black!10} 156 & 75 & 832 & 3457 & 798$^*$ & 5.10 & 0.63 & 1.45 & 0.00 & 609.13 & 103.84 & 0.00 \\ 
 157 & 75 & 832 & 3457 & 821$^*$ & 5.87 & 1.34 & 1.32 & 0.00 & 301.01 & 39.39 & 0.00 \\ 
\rowcolor{black!10} 158 & 75 & 832 & 3457 & 816$^*$ & 5.66 & 0.49 & 1.54 & 0.00 & 522.14 & 28.66 & 0.00 \\ 
 159 & 75 & 832 & 3457 & 820$^*$ & 5.12 & 0.24 & 81.59 & 0.12 & 654.14 & 49.56 & 0.12 \\ 
\rowcolor{black!10} 160 & 75 & 832 & 3457 & 815$^*$ & 5.91 & 1.47 & 68.23 & 0.12 & 297.86 & 30.99 & 0.00 \\ 
 161 & 75 & 832 & 13828 & 903 & 14.42 & 25.25 & 1977.61 & -0.55 & 1071.69 & 211.42 & 0.44 \\ 
\rowcolor{black!10} 162 & 75 & 832 & 13828 & 953 & 13.11 & 18.05 & 3716.49 & 0.63 & 943.99 & 136.26 & -2.52 \\ 
 163 & 75 & 832 & 13828 & 892 & 13.52 & 16.03 & 3795.01 & 2.02 & 1221.56 & 260.01 & 0.56 \\ 
\rowcolor{black!10} 164 & 75 & 832 & 13828 & 915 & 14.10 & 24.59 & 3343.98 & 1.75 & 1013.60 & 153.48 & 0.98 \\ 
 165 & 75 & 832 & 13828 & 896 & 13.90 & 17.97 & 2741.60 & 4.24 & 434.95 & 61.89 & 1.12 \\ 
\rowcolor{black!10} 171 & 75 & 1110 & 6155 & 787$^*$ & 12.41 & 0.13 & 9.65 & 0.00 & 4.66 & 2.63 & 0.00 \\ 
 172 & 75 & 1110 & 6155 & 785$^*$ & 10.83 & 0.13 & 50.48 & 0.00 & 257.80 & 46.85 & 0.00 \\ 
\rowcolor{black!10} 173 & 75 & 1110 & 6155 & 783$^*$ & 12.67 & 2.17 & 2.18 & 0.00 & 1297.35 & 118.07 & 0.00 \\ 
 174 & 75 & 1110 & 6155 & 784$^*$ & 12.30 & 0.64 & 3.04 & 0.13 & 391.55 & 43.63 & 0.00 \\ 
\rowcolor{black!10} 175 & 75 & 1110 & 6155 & 797$^*$ & 11.92 & 1.51 & 3.41 & 0.00 & 33.91 & 7.69 & 0.00 \\ 
 176 & 75 & 1110 & 24620 & 867 & 27.07 & 19.49 & 2841.94 & -0.35 & 1290.64 & 240.51 & 1.15 \\ 
\rowcolor{black!10} 177 & 75 & 1110 & 24620 & 851 & 29.82 & 23.15 & 3817.10 & 3.64 & 1318.79 & 85.64 & 2.82 \\ 
 178 & 75 & 1110 & 24620 & 892 & 30.08 & 23.09 & 2969.99 & 6.84 & 1115.62 & 162.08 & -1.23 \\ 
\rowcolor{black!10} 179 & 75 & 1110 & 24620 & 864 & 30.47 & 24.54 & 3356.33 & 2.20 & 1373.69 & 150.56 & 11.81 \\ 
 180 & 75 & 1110 & 24620 & 882 & 29.08 & 21.66 & 3096.31 & -0.45 & 756.48 & 32.41 & 3.74 \\ 
\rowcolor{black!10} 183 & 75 & 1110 & 43085 & 1194 & 34.48 & 100.00 & 3619.07 & 100.00 & 1699.45 & 47.18 & 100.00 \\ 
 186 & 100 & 990 & 4896 & 1119$^*$ & 11.52 & 3.93 & 141.09 & 0.00 & 834.38 & 48.88 & 0.00 \\ 
\rowcolor{black!10} 187 & 100 & 990 & 4896 & 1137$^*$ & 10.44 & 1.67 & 168.84 & 0.00 & 408.35 & 63.67 & 0.00 \\ 
 188 & 100 & 990 & 4896 & 1113$^*$ & 10.59 & 2.70 & 138.40 & 0.00 & 1537.76 & 266.18 & 0.09 \\ 
\rowcolor{black!10} 189 & 100 & 990 & 4896 & 1110$^*$ & 11.09 & 4.05 & 104.98 & 0.09 & 684.33 & 141.51 & 0.00 \\ 
 190 & 100 & 990 & 4896 & 1090$^*$ & 10.54 & 2.20 & 11.84 & 0.09 & 1284.48 & 167.24 & 0.09 \\ 
\rowcolor{black!10} 191 & 100 & 990 & 19583 & - & 23.31 & 100.00 & 3847.28 & 100.00 & 1290.60 & 140.05 & -100.00 \\ 
 192 & 100 & 990 & 19583 & 1491 & 24.13 & 100.00 & 3832.03 & 100.00 & 1159.85 & 44.61 & 17.91 \\ 
\rowcolor{black!10} 193 & 100 & 990 & 19583 & 1510 & 23.93 & 100.00 & 3850.73 & 100.00 & 1151.45 & 43.71 & 14.97 \\ 
 194 & 100 & 990 & 19583 & 1441 & 23.94 & 100.00 & 3980.25 & 100.00 & 1149.78 & 42.00 & 23.73 \\ 
\rowcolor{black!10} 195 & 100 & 990 & 19583 & 1560 & 24.58 & 100.00 & 3939.26 & 100.00 & 1344.19 & 174.61 & 9.49 \\ 
 201 & 100 & 1485 & 11019 & 1079$^*$ & 28.92 & 5.28 & 150.34 & 0.00 & 1898.20 & 301.64 & 0.19 \\ 
\rowcolor{black!10} 202 & 100 & 1485 & 11019 & 1056$^*$ & 32.30 & 2.65 & 107.84 & 0.00 & 935.21 & 183.48 & 0.00 \\ 
 203 & 100 & 1485 & 11019 & 1059$^*$ & 29.54 & 1.70 & 8.75 & 0.28 & 761.93 & 100.02 & 0.00 \\ 
\rowcolor{black!10} 204 & 100 & 1485 & 11019 & 1046$^*$ & 30.32 & 1.24 & 115.85 & 1.05 & 871.40 & 119.46 & 0.00 \\ 
 205 & 100 & 1485 & 11019 & 1072$^*$ & 28.53 & 3.45 & 261.44 & 0.00 & 1258.39 & 43.68 & 0.00 \\ 
\rowcolor{black!10} 206 & 100 & 1485 & 44075 & 1374 & 70.00 & 100.00 & 3846.93 & 100.00 & 1252.91 & 42.30 & 10.04 \\ 
 207 & 100 & 1485 & 44075 & 1291 & 66.76 & 100.00 & 3848.22 & 100.00 & 1173.60 & 30.59 & 21.84 \\ 
\rowcolor{black!10} 208 & 100 & 1485 & 44075 & 1344 & 66.37 & 100.00 & 3836.87 & 100.00 & 1277.02 & 34.77 & 24.55 \\ 
 209 & 100 & 1485 & 44075 & 1286 & 66.20 & 100.00 & 3858.32 & 100.00 & 1234.78 & 45.74 & 20.37 \\ 
\rowcolor{black!10} 210 & 100 & 1485 & 44075 & 1370 & 62.92 & 100.00 & 3844.58 & 100.00 & 1162.74 & 30.23 & 14.60 \\ 
 216 & 100 & 1980 & 19593 & 1031$^*$ & 65.49 & 3.69 & 140.15 & 0.00 & 1478.66 & 16.70 & 0.00 \\ 
\rowcolor{black!10} 217 & 100 & 1980 & 19593 & 1036$^*$ & 70.46 & 4.34 & 95.81 & 0.00 & 275.72 & 46.57 & 0.00 \\ 
 218 & 100 & 1980 & 19593 & 1024$^*$ & 66.73 & 3.91 & 122.80 & 0.00 & 761.29 & 101.41 & 0.00 \\ 
\rowcolor{black!10} 219 & 100 & 1980 & 19593 & 1025$^*$ & 61.12 & 1.76 & 82.52 & 0.00 & 1795.81 & 294.19 & 0.00 \\ 
 220 & 100 & 1980 & 19593 & 1028$^*$ & 68.04 & 4.67 & 153.57 & 0.00 & 1064.45 & 159.70 & 0.00 \\ 
\rowcolor{black!10} 221 & 100 & 1980 & 78369 & 1234 & 140.14 & 100.00 & 3893.96 & 100.00 & 1516.22 & 102.96 & 20.10 \\ 
 222 & 100 & 1980 & 78369 & 1187 & 149.79 & 41.53 & 3697.46 & 100.00 & 1283.16 & 54.50 & 25.61 \\ 
\rowcolor{black!10} 223 & 100 & 1980 & 78369 & 1213 & 139.31 & 100.00 & 3999.99 & 100.00 & 1310.06 & 90.75 & 21.85 \\ 
 224 & 100 & 1980 & 78369 & 1221 & 140.17 & 100.00 & 3999.16 & 100.00 & 1311.85 & 45.53 & 17.85 \\ 
\rowcolor{black!10} 225 & 100 & 1980 & 78369 & 1245$^*$ & 127.94 & 35.42 & 3856.70 & 100.00 & 1274.09 & 43.89 & 11.65 \\ \hline
\rowcolor[HTML]{C0C0C0} \multicolumn{5}{!{\color{black!20}\vrule width 0.8pt}c!{\color{black!20}\vrule width 0.8pt}}{Average values} &  \textbf{14.16} & 23.05 & 924.04 & 12.72 & 440.83 & 47.26 & \textbf{1.73} \\ 
\rowcolor[HTML]{C0C0C0} \multicolumn{5}{!{\color{black!20}\vrule width 0.8pt}c!{\color{black!20}\vrule width 0.8pt}}{\%best} & \multicolumn{2}{c!{\color{black!20}\vrule width 0.8pt}}{29.80} & \multicolumn{2}{c!{\color{black!20}\vrule width 0.8pt}}{49.01} & \multicolumn{3}{c!{\color{black!20}\vrule width 0.8pt}}{\textbf{74.83}}\\
\multicolumn{11}{c}{\label{tab:kernel_cmst}}
\end{longtable}}\vspace*{-5mm}
It is clear from the comparison between the results of our KS approach and the classic KS one that using the maximum independent set problem, as well as the {heuristic} algorithm, is crucial to address the MSTC problem. The increase in the computational time with respect to other existing methods is outweighed by the considerable reduction in the performance gap. 

\end{document}